\documentclass[10pt,twosided]{article}
\usepackage[tbtags]{amsmath}
\usepackage{amssymb}
\allowdisplaybreaks[4]
\usepackage{amsthm}

\usepackage{color}

\usepackage{amssymb,color}

\definecolor{c20}{rgb}{0.,0.7,0.}
\definecolor{c30}{rgb}{0.,0.,0.}
\definecolor{c40}{rgb}{1,0.1,0.7}
\definecolor{c50}{rgb}{1,0,0}
\definecolor{c60}{rgb}{0,0.9,0.1}

\def\eeH#1{\textcolor{c50}{#1}}
\def\eeH#1{#1}
\def\aH#1{\textcolor{c50}{#1}}
\def\cH#1{#1}

\def\aH#1{#1}

\usepackage{amssymb,color}

\newcommand{\kb}[1]{\boldsymbol{#1}}
\newcommand{\vk}[1]{\kb{#1}}

\newcommand{\abs}[1]{\left\lvert #1 \right\rvert}

\newcommand{\E}[1]{\mathbb{E}\left(#1\right)}

\newcommand{\pk}[1]{\mathbb{P} \left( #1 \right) }
\newcommand{\pc}[1]{\mathbb{P}\left(#1 \right)}

\newcommand{\R}{\mathbb{R}}

\newcommand{\todis}{\stackrel{d}{\to}}

\newcommand{\BQN}{\begin{eqnarray}}
\newcommand{\EQN}{\end{eqnarray}}
\newcommand{\BQNY}{\begin{eqnarray*}}
\newcommand{\EQNY}{\end{eqnarray*}}

\newcommand{\BS}{\begin{sat}}
\newcommand{\ES}{\end{sat}}
\newcommand{\BT}{\begin{theo}}
\newcommand{\ET}{\end{theo}}
\newcommand{\BK}{\begin{korr}}
\newcommand{\EK}{\end{korr}}

\newcommand{\BD}{\begin{de}}
\newcommand{\ED}{\end{de}}

\newcommand{\BIT}{\begin{itemize}}
\newcommand{\EIT}{\end{itemize}}
\newcommand{\BDI}{\begin{description}}
\newcommand{\EDI}{\end{description}}

\newcommand{\BRMS}{\begin{remarks}}
\newcommand{\ERMS}{\end{remarks}}
\newcommand{\BRM}{\begin{remark}}
\newcommand{\ERM}{\end{remark}}
\newcommand{\BEL}{\begin{lem}}
\newcommand{\EEL}{\end{lem}}

\newtheorem{theo}{Theorem}[section]
\newtheorem{sat}[theo]{Proposition}
\newtheorem{de}[theo]{Definition}
\newtheorem{lem}[theo]{Lemma}

\newtheorem{korr}[theo]{Corollary}
\newtheorem{remark}[theo]{Remark}
\newtheorem{remarks}[theo]{Remarks}

\newcommand{\nelem}[1]{{Lemma \ref{#1}}}

\newcommand{\netheo}[1]{{Theorem \ref{#1}}}

\newcommand{\prooftheo}[1]{ \textbf{Proof of Theorem} \ref{#1}: }

\newcommand{\proofkorr}[1]{\textbf{Proof of Corollary} \ref{#1}:}

\newcommand{\COM}[1]{}

\newcommand{\QED}{\hfill $\Box$}

\def\invf{\overleftarrow{f}}
\def\invfs{\overleftarrow{f^*}}
\def\invF{\overleftarrow{F}}

\def\rw{\rightarrow}

\def\IF{\infty}

\date{}

\def\LT{\left}
\def\RT{\right}

\def\Cov{\mathrm{Cov}}

\def\Corr{\mathrm{Corr}}

\def\NN{\mathbb{N}}
\def\finv{\overleftarrow{f}}

\def\EH#1{\textcolor{c30}{#1}}
\def\EH#1{#1}
\def\HH#1{\textcolor{c20}{#1}}
\def\HH#1{#1}
\begin{document}

\title{\bf Extremes of Locally Stationary Chi-square Processes with Trend}

\bigskip
\author{Peng Liu\thanks{Department of Actuarial Science, University of Lausanne, UNIL-Dorigny 1015 Lausanne, Switzerland}
\thanks{Mathematical Institute, University of Wroc\l aw, pl. Grunwaldzki 2/4, 50-384 Wroc\l aw, Poland}\quad and Lanpeng Ji$^*$\thanks{Institute for Information and Communication Technologies, School of Business and Engineering Vaud (HEIG-VD), University of Applied Sciences of Western Switzerland, Switzerland}
}

 \maketitle
\vskip -0.61 cm

  \centerline{\today{}}

\bigskip
{\bf Abstract:}
Chi-square processes with trend appear naturally as limiting processes in various statistical models. In this paper
we are concerned with %cases for which the supremum taken over $(0,1)$ of such processes is finite only for particular admissible trends. Our main findings are about
the exact tail asymptotics of the supremum taken over $(0,1)$ of a class of locally stationary chi-square processes with particular admissible trends.
An important tool for establishing our results is a weak version of Slepian's lemma  for chi-square processes. %As a by product, w
%We also show that the admissibility condition  on the trend function  given in terms of the local structure of the involved Gaussian process %. For some special cases, such as Brownian bridge and Bessel process, the relation between the admissibility condition and is closely related to the Kolmogorov-Dvoretsky-Erd\H{o}s integral test of the chi-square process. 
Some special cases including squared Brownian bridge and Bessel process are discussed.



%%%%%%%%%%%%%%%%%%%%%%%%%%%%%%%%%%%5555555
\section{Introduction}

Let $\widehat{G}_n(t), t\in (0,1)$ be  the empirical distribution of $n$ independent random variables with uniform distribution on $(0,1)$ and define the following test statistic
\BQNY
L_{n,E}^\nu:= \sup_{t\in E}\left(n{K}(\widehat{G}_n(t),t)-g_\nu(t)\right),\ \ E=(0,1), %{K}_n(s,t)=n\Bigl(s \ln\frac{s}{t}+(1-s)\ln\frac{1-s}{1-t}\Bigr), \quad E:=(0,1),
\EQNY
where  ${K}(s,t)=s \ln\frac{s}{t}+(1-s)\ln\frac{1-s}{1-t} $,  and
$g_\nu(t), t>0,\nu>0$ is a trend function defined by
$$g_\nu(t)={c}(t)+\nu\ln(1+{c}^2(t)), \quad {c}(t)=\ln \bigl(1- \ln(4t(1-t))\bigr), \quad \quad t\in E.$$
%, $t\in(0,1)$.\\
Referring to Theorem  3.2 in \cite{Lutz}  we have for any $\nu>\aH{3/4}$ %(\aH{we conjecture that for $\nu>3/4$, the following argument is also valid but this wasn't proved in  \cite{Lutz}}) 
%\aH{the following convergence  in distribution holds:}
 \BQN\label{INBR}
2L_{n,E}^\nu &\todis & L_{E}^\nu:=\sup_{t\in E}\Bigl( \chi^2(t)-2g_\nu(t)\Bigr),\quad n\to \IF
 \EQN
holds, where $\chi (t)= B(t)/ \sqrt{t(1-t)}$ is the normalized standard Brownian bridge. Furthermore, as shown in Theorem 3.4 in \cite{Lutz}, the convergence to $L_E^\nu$ also  holds for another important test statistic. Note that %the condition $\nu>1$ is needed since $E=(0,1)$. Actually,
 as discussed therein $L_{E}^\nu<\IF$ almost surely (a.s.) for any $\nu>3/4$.
 % Clearly, when $E \subset (0,1)$ is \cH{a compact interval}, then $L_E^\nu< \IF$ a.s., which can be easily shown  for any continuous trend function $g_\nu(\cdot)$. % can be chosen to be quite general.\\

The supremum type test statistics appear naturally  in different contexts in  statistical problems; see also  the recent contributions  \cite{JA1,BB,JA2, JAPIT}. %In the aforementioned contributions, it is shown that "trimmed versions" (supremum is taken on some compact interval of $(0,1)$) of the supremum of chi-square processes (without trend) is a candidate for some limiting statistics, this is due to fact that the "non-trimmed versions" are usually infinite. Motivated by
%Based on the aforementioned contributions, the limiting process  $\chi^2$ can be more general, for instance, it can be assumed to be a chi-square process. Since, even for the case that $X$ is a Brownian bridge it is not possible to derive some explicit formula for the distribution function of the random variable $L_E$,  of interest is to derive its tail behaviour.
  In statistical applications, of interest is to determine the critical values of the test statistics, which is usually difficult since the distribution of the statistics is unknown. It is thus important to obtain approximation of them based on the tail asymptotic behavior of the limit  of the statistic. For example, in our context it is of interest to know  the asymptotic behavior of $\pk{L_E^\nu>u}$ as $u\to\IF$ for any $\nu>\aH{3/4}$. % (or if possible $\nu>3/4$).

Note that $L_E^\nu$ is the supremum of a simple chi-square process (with 1-degree  of freedom) with trend. Numerous contributions have been  devoted to the study of the tail asymptotics of the supremum of chi or chi-square processes; see e.g., \cite{HAJI2014, konstantinides2004gnedenko,Lindgren1989,chiLiu,Pit96,tanH2012} and the references therein. So far in the literature there are no results available which can be used for the derivation of the tail behavior of $L_E^\nu$ since it is a supremum taken over an open interval  of a (non-Gaussian) chi-square process with trend. Given this significant gap and the fact that chi-square process with trend appears naturally in numerous statistical problems and applications, in this paper we shall focus on the tail asymptotics of the supremum (taken over an open interval) of a large class of chi-square processes with trend. More precisely, define
\BQN\label{INCHI}
 \chi_{\vk{b}}^2(t)= \sum_{i=1}^n b_i^2 X_i^2(t), \quad t\in (0, \IF),
 \EQN
  where
  $$1=b_1= \cdots= b_k> b_{k+1} \ge \cdots \ge b_n>0$$
   and $X_i$'s are independent copies of a centered
  Gaussian process $X$ with a.s. continuous sample paths. We are interested in the asymptotics of
  \BQN \label{chig}
  \pk{\sup_{t\in E}\Bigl(\chi_{\vk{b}}^2(t)-g(t)\Bigr)>u} %,\ \ E=(0,1)
  \EQN
as $u\to\IF$, for certain Gaussian processes $X$ and nonnegative continuous trend functions $g(\cdot)$. \aH{Restrictions on $X$ and $g(\cdot)$  will be specified to first ensure  that} %the existence of the asymptotics, which in turn implies that
\BQN \label{chigfinite}
\sup_{t\in E}\Bigl(\chi_{\vk{b}}^2(t)-g(t)\Bigr)<\IF,\ \ a.s.
\EQN
holds.

In the special case that  $n=1$ and $X$ is the normalized standard Brownian bridge, then $\chi_1=X$, the same as $\chi$ appearing in \eqref{INBR}, is a locally stationary Gaussian process in the notion of \cite{Berman74, Husler90}.
Precisely,   %for this special case note that %$\overline{X}(t)=\frac{X(t)}{\sqrt{t(1-t)}}$ satisfies,
%for any $t \in (0,1)$
\BQNY%\label{UC}
\lim_{h\to 0}\frac{1-\E{ \chi(t)\chi(t+h)}}{\abs{h}}=C(t)
\EQNY
holds uniformly in $t\in I$, any compact interval in $(0,1)$, where $C(t)= \frac{1}{2t(1-t)}$ satisfying $C(0)=C(1)=\IF$.
%Moreover, the convergence in \eqref{UC}  is uniform on any compact interval in $(0,1)$.
Motivated by this fact in this paper we shall consider a large class of centered locally stationary Gaussian processes $\{X(t),t\in (0,1)\}$ with a.s. continuous sample paths,  unit variance and correlation function $r(\cdot,\cdot)$ such that
\BQN\label{INCO}
\lim_{h\to 0}\frac{1-r(t,t+h)}{K^2(\abs{h})}=C(t)
\EQN
uniformly in $t\in I$, any compact interval in $(0,1)$, where $K(\cdot) $ is a positive regularly varying function at 0 with index $\alpha/2\in(0, 1]$, and $C(\cdot)$ is a positive continuous function satisfying $C(0)=\IF$ or $C(1)=\IF$.

It is noted that condition \eqref{INCO} for the correlation function of $X$ seems natural and is satisfied by many other interesting Gaussian processes. For instance,
let $\{B_H(t), t\ge0\}$  be the standard fractional Brownian motion (fBm) with Hurst index $H\in(0,1)$ and covariance function given by
$$
\Cov(B_H(s),B_H(t))=\frac{1}{2}(t^{2H}+s^{2H}-\abs{t-s}^{2H}),\ \ \ s,t\ge0.
$$
We can show that the normalized fBm $\{B_H(t)/t^{H},t\in(0,\IF)\}$ is a locally stationary Gaussian process. In fact, the correlation of the normalized fBm satisfies
\BQN \label{eq:corr_fBm}
\lim_{h\to0}\frac{1-\text{Corr}\LT(\frac{B_H(t)}{t^{H}},\frac{B_H(t+h)}{(t+h)^{H}}\RT)}{|h|^{2H}}= \frac{1}{2t^{2H}} %, \ \ h\rw 0,
\EQN
uniformly in $t$, for any compact interval in $(0,\IF)$. %, where $C(t)= \frac{1}{2t^{2H}}$ satisfying $C(0)=\IF, C(1)=1/2$. %Note that in this case $1$ is called a {\it singular point} of $C(\cdot)$.

%%%%%%%%%%%%%%%%%%%%%%%%%%%%%%%%%%%%%%%%%%%%%%%%%%%%%%%%%%%%%%%%%%%%%%%%%%%%%%%%%%%%%%%%%%%%%%%%%%%%%%%%%%%%%%%%%%%%%%%%%%%%%%%%%%%55555
\COM{
let $\{B_{F,H}(t),t\ge0\},F,H\in (0,1]$ be a bi-fractional Brownian motion (bi-fBm) with covariance function
\BQNY
\Cov\left(B_{F,H}(t), B_{F,H}(s)\right)=\frac{1}{2^{F}}\left((t^{2H}+s^{2H})^{F}-|t-s|^{2FH}\right), ~~t,s\geq 0.
\EQNY
%Note that when $K_0=1$, $\{B_{K_0,H}(t), t\geq 0\}$ is a fractional Brownian motion with index $H$.
The standardized bi-fBm process $X_{F,H}(t)=t^{-FH}B_{F,H}(t)$, $t_0 \in (0,T)$ satisfies
\BQNY
\mathbb{E}\left(X_{F,H}(t)X_{F,H}(s)\right)=1-\frac{1}{2^{F}t_0^{2FH}}|t-s|^{2FH}(1+o(1)), ~~t,s\rw t_0\in(0,T),
\EQNY
and hence condition \eqref{INCO} holds.
%The coefficient $\frac{1}{2^{K_0}t_0^{2K_0H}}$ tends to infinity as $t_0$ tends to zero.\\
Further, let  $\{W(t), t\geq 0\}$ be a centered stationary Gaussian process with continuous sample paths, unit variance
and correlation function $r(\cdot)$ satisfying $r(t)=1-K^2(t)(1+o(1)), t\rw 0$ with $K^2(\cdot)$ a positive regularly varying function at 0 with index $\alpha \in (0,2]$. Then $X(t)= W(\ln t)$ satisfies the assumptions of \netheo{Th0} since \eqref{INCO} holds with $C(t)= t^{- \alpha}$, namely, for any $t_0\in (0,T)$
%.A larger class of Guassian processes that is obtained by s discussed in [], we consider a Gaussian process $\{X(t)=W(\ln t), t>0\}$ $
\BQNY
\mathbb{E}\left(X(t)X(s)\right)=1-\frac{1}{t_0^\alpha}K^2(|t-s|),~~t,s\rw t_0\in(0,T).
\EQNY
}
%%%%%%%%%%%%%%%%%%%%%%%%%%%%%%%%%%%%%%%%%%%%%%%%%%%%%%%%%%%%%%%%%%%%%%%%%%%%%%%%%%%%%%%%%%%%%555555

In order to derive the asymptotics of \eqref{chig} for $E=(0,1)$,  we need to impose several other conditions on $g(\cdot)$ and on the locally stationary Gaussian process $X$; see Section 2. One natural restriction on $g(\cdot)$ is to ensure that it satisfies \eqref{chigfinite}; functions $g(\cdot)$ satisfying \eqref{chigfinite} are thus called {\it admissible functions}. \aH{ \netheo{Addition}  below shows that the admissibility of functions $g(\cdot)$ is related to a generalized Kolmogorov-Dvoretsky-Erd\H{o}s integral test (or the law of iterated logarithm) of the corresponding processes; see also Appendix.  For instance, as discussed in Corollary \ref{cor} below $2g_v(\cdot)$ in \eqref{INBR} is admissible if and only if $\nu>3/4$; see also \cite{Lutz}.} In addition, as an application of our main result given in \netheo{TH1}, we obtain the asymptotics of $\pk{L_E^\nu>u}$ as $u\to\IF$ for any $\nu>3/4$. Other  examples related to the fBm and the Bessel process are also discussed.

%As shown in  \netheo{TH1}  below, the tail asymptotics of $L_E$ can be recovered for several interesting cases of $X$ and trends $g$; special interesting examples are $X$ being a bi-fractional Brownian motion,
 %sub-fractional Brownian motion, or  mean integrated fractional Brownian motion. \\
Organization of the rest of the paper: In Section 2 we present our main results which are illustrated by several examples. Further results are discussed in Section 3 for the case where the Gaussian processes $X_i$'s are not identically distributed.  All the proofs are relegated to Section 4, whereas some  technical results are postponed to Appendix.

\section{Main Results}
We begin with some preliminary notation. We shall use the standard   notation for asymptotic equivalence of two functions $f(\cdot)$ and $h(\cdot)$. That is, for any $x_0\in\R\cup\{\IF\},$  write
$ f(x)=h(x)(1+o(1))$  or simply $f(x)\sim h(x)$, if  $ \lim_{x \to x_0}  {f(x)}/{h(x)} = 1$, and write $ f(x) = o(h(x))$, if $ \lim_{x \to x_0}  {f(x)}/{h(x)} = 0$.  Denote by $\Gamma(\cdot)$
the Euler Gamma function. Further, denote by $\overleftarrow{K}(\cdot)$  the generalized inverse function of $K(\cdot)$, and set $q(u)= \overleftarrow{K}(u^{-1/2})$ for any $u>0$. Additionally, write
 $\mathcal{H}_\alpha,\EH{\alpha \in (0,2]}$  for the  Pickands constant; see \cite{DHE16, debicki2008note, debicki2002ruin,
 nonhomoANN, dekos14, DM15, DikerY, PicandsA,  Pit96} for its definition and generalizations.

Our first result is concerned with the asymptotics of \eqref{chig} for the case that $E$ is any compact sub-interval in  $(0,1)$. The trend function appears in the
asymptotics indirectly through the following constant
\BQN \label{Jg}
J_E^g:= \int_{t \in E}C^{1/\alpha}(t)e^{-\frac{g(t)}{2}}dt<\IF.
\EQN

\BT\label{Th0}
 Let $\{X(t),t\in E\}$, with $E$ a compact interval in $(0,1)$, be  a centered locally stationary Gaussian process  with a.s. continuous sample paths, unit variance  and
correlation function $r(\cdot,\cdot)$ that satisfies \eqref{INCO} and $r(s,t)<1$, $s\neq t \in E$. Then, for any nonnegative continuous trend function $g(\cdot)$
we have
\BQN\label{PRC}
\pk{\sup_{t\in E}\Bigl(\chi_{\vk{b}}^2(t)-g(t)\Bigr)>u}\sim\mathcal{H}_{\alpha}G_{\vk{b}}J_E^g \eeH{\frac{u^{k/2-1}e^{-u/2}}{q(u)}}
\EQN
as $u\to \IF$,
\eeH{where  $G_{\vk{b}}=\frac{2^{1-k/2}}{\Gamma(k/2)}\prod_{i=k+1}^n\left(1-b_i^2\right)^{-1/2}$ and \aH{$q(u)= \overleftarrow{K}(u^{-1/2})$} (with the convention $\prod_{i=l}^m=1$ if $m<l$)}.
\ET

 \begin{remark}\label{remCI}
a) We see from \netheo{Th0} that if $X$ is stationary and $g(t)\equiv 0$ the above result coincides with that derived in \cite{Piterbarg94}. In \cite{PitSta01} the authors obtained the tail asymptotics of the supremum of a class of Gaussian random fields with trend indexed on smooth manifolds. It is worth noting that in \netheo{Th0} above the trend function $g(\cdot)$ contributes to the asymptotics through $J_E^g$, which is \aH{quite} different from that in the aforementioned paper. \\
b) Clearly, if $X$ is well-defined on $(0,1]$ (e.g., the normalized fBm, see \eqref{eq:corr_fBm}), then under the assumptions of Theorem \ref{Th0} we have \eqref{PRC} holds for any  $E=[a,1]$  with $a\in(0,1)$.
In fact, the set $E$ in the above theorem can be any compact interval in $\R$ provided that everything is well-defined on it.
 \end{remark}
%In our notion of the locally stationary Gaussian process, we assume  that $C(0)=\IF$ or $C(1)=\IF$. In such a case, one may wonder if  $\int_{0}^1 (C(t))^{1/\alpha} dt<\IF$ implies \eqref{PRC} holds
Our main \aH{targets below are} to find out under which conditions \aH{(\ref{chigfinite}) holds and in such a case to establish \eqref{PRC}} by passing from the compact interval $E$  to $(0,1)$. Of course, \aH{ conditions should be imposed} on the local structures of the process $X$ and $g(\cdot)$ at 0 and 1, separately. %Note that as discussed in Remark \ref{remCI} if $X$ is well defined on $(0,1]$ we only need to consider the conditions around 0. %Since we aim to focus on the case where $E=(0,1)$, we shall implicitly assume that the involved process $X$ is defined only on $(0,1)$. Our results can also be applied to the case with $E=(0,1]$ as illustrated in ...
\aH{According to   the proof of Theorem \ref{THM} in the Appendix (see (\ref{con1}) and (\ref{con2})), we believe that, under some mild conditions, % for correlation of $X$ and trend function $g$, 
the sufficient and necessary condition for (\ref{chigfinite}) to hold is  given as 
\BQN\label{cs}
I_g(S):=\left|\int_{1/2}^S (C(t))^{1/\alpha}\frac{(g(t))^{\frac{k}{2}-1}}{q(g(t)) }e^{-\frac{g(t)}{2}}dt\right|<\IF, \quad S\in \{0,1\}.
\EQN
Furthermore,  under these conditions we  show that \eqref{PRC} holds with $E=(0,1)$. 
Moreover, note that if $\left|\int_{1/2}^S (C(t))^{1/\alpha}dt\right|<\IF$, then (\ref{cs}) holds for all $g(\cdot)$ satisfying $\lim_{t\to S} g(t)=\IF$. % implying that we don't need the complicated assumption (\ref{cs}). Besides, 
It turns out that  in this simpler case %the assumptions on correlation of $X$ and the proof are both much simpler than and different from the case $\left|\int_{1/2}^S (C(t))^{1/\alpha}dt\right|=\IF$.  
the conditions imposed on the correlation function of $X$  and the proof of main results can be significantly simplified. Thus, to simplify argumentation,}
\COM{First, let us give some intuition based on the asymptotics obtained in \eqref{PRC}.
 One could conjecture that if  $\int_{0}^1 (C(t))^{1/\alpha} dt<\IF$ (implying $J_{(0,1)}^g<\IF$), then we should have \eqref{PRC}
for $E=(0,1)$. However, \aH{we believe} that this simple condition solely is not sufficient for the derivation of \eqref{PRC} with $E=(0,1)$, and additional condition is needed; see condition {\bf D} below.  %under which  we shall assume that $J_{(0,1)}^g<\IF,$,
Further, if $\int_{0}^1 (C(t))^{1/\alpha} dt=\IF$, then $J_{(0,1)}^g<\IF$ is also not  sufficient for the validity of \eqref{PRC} with $E=(0,1)$. In this case a stronger condition (see condition {\bf C} below)  is necessary, which in turn gives a nice criterion for admissible trends $g(\cdot)$ (see \netheo{Addition}).}
 different scenarios will be discussed according to whether $\int_{0}^{1/2}(C(s))^{1/\alpha}ds<\IF$ and whether $\int_{1/2}^1(C(s))^{1/\alpha}ds=\IF$, and different additional assumptions are needed accordingly.
 For this purpose of crucial importance is the following function
\BQNY
f(t)=\int_{\cH{1}/2}^t(C(s))^{1/\alpha}ds,\ \ \ t\in(0,1).
\EQNY
We denote by $\invf(t), t\in(f(0), f(1))$ the inverse function of $f(t),t\in(0,1)$. Further, for any $d>0$,
 let $s_{j,d}^{(1)}=\invf(jd)$, $j\in \mathbb{N}\cup \{0\}$ if $f(1)=\IF$, and let $s_{j,d}^{(0)}=\invf(-jd)$, $j\in \mathbb{N}\cup\{0\}$ if $f(0)=-\IF$. Denote $\Delta_{j,d}^{(1)}=[s_{j-1,d}^{(1)},s_{j,d}^{(1)}], j\in \mathbb{N}$ and
  $\Delta_{j,d}^{(0)}=[s_{j,d}^{(0)},s_{j-1,d}^{(0)}], j\in \mathbb{N}$, which give a partition of $[1/2,1)$ in the case $f(1)=\IF$ and a partition of $(0,1/2]$ in the case $f(0)=-\IF$, respectively.
 % then we can derive a division of $[T/2,T)$ denoted by $\Delta_{j,d}^{(T)}=[s_{j-1,d}^{(T)},s_{j,d}^{(T)}], j\in \mathbb{N}$. Similarly, if, we can derive a division for, i.e., $\Delta_{j,d}^{(0)}=[s_{j,d}^{(0)},s_{j-1,d}^{(0)}], j\in \mathbb{N}$ with.
 % We first consider the case that
%$$\cH{I(E):=\int_{t\in E}}(C(t))^{1/\alpha}dt=\IF.$$

In addition to the local stationarity of $X$ in \eqref{INCO}, we need to impose the following (scenario-dependent) restrictions on
the trend $g(\cdot)$ and the correlation function $r(\cdot,\cdot)$. Let therefore $S\in\{0,1\}$.  \\
\textbf{Condition A}($S$):  The  trend function $g(\cdot)$  is monotone in a neighborhood of $S$ and  $\lim_{t\to S} g(t)=\IF$.\\
%\EH{$\lim_{t\rw \cH{S} }g(t)=\IF$} and there exist a constant $0<\delta_g<T/2$  such that $g$ is monotone over $E\cap (S-\delta_g, S+\delta_g)$ \cH{for $S\in \{0,T\}$}. \\
 \textbf{Condition B}$(S)$: Suppose that there exists some constant  $d_0>0$ such that
\BQNY
  \limsup_{j\rw \IF}\sup_{t\neq s\in\Delta^{(S)}_{j,d_0}}\frac{1-r(t,s)}{K^2(|f(t)-f(s)|)}\cH{ < \IF},
 \EQNY
and when $\alpha=2$ and $k=1$, assume further
 \BQN\label{K}
 K^2(|t|)=O(t^2), \ \ t\rw 0.
 \EQN
\textbf{Condition C}$(S)$:  \aH{With $I_g(S)$ defined in \eqref{cs} it holds that}
\BQN\label{INAS} 
I_g(S)%=\left|\int_{1/2}^S (C(t))^{1/\alpha}\frac{(g(t))^{\frac{k}{2}-1}}{q(g(t)) }e^{-\frac{g(t)}{2}}dt\right|
<\IF.
 %~~ \int_{0}^{d\delta_0}\widetilde{R}_{d}(t)(1+g(t))^{(\frac{k}{2}-1)_{+} +1+\eta_1}e^{-\frac{g(T-t)}{2}}dt<\IF.
 \EQN
 \textbf{Condition D}($S$): The following is satisfied  % For any $S\in\{0,T\}$   %There exists a constant $M_F$ such that
\BQNY
\limsup_{\delta\rw 0}\sup_{t\neq s\in (0,\delta)}\frac{1-r(|S-t|,|S-s|)}{K^2(|f(|S-t|)-f(|S-s|)|)} \cH{< \IF}.
\EQNY
\textbf{Condition E}$(S)$: Suppose that there exists some constant  $d_0>0$ such that
\BQNY
    \liminf_{j\rw \IF}\inf_{t\neq s\in\Delta^{(S)}_{j,d_0}}\frac{1-r(t,s)}{K^2(|f(t)-f(s)|)}> 0.
 \EQNY
 Moreover, there exist  $j_0, l_0\in \mathbb{N}$, $M_0, \beta>0, $ such that for $j\geq j_0,$ $l\geq l_0$,
 \BQN\label{B2}
 \sup_{s\in \Delta^{(S)}_{j+l,d_0},t\in \Delta^{(S)}_{j,d_0} }|r(s,t)|<M_0l^{-\beta}.
 \EQN
Let
$$
\mathcal{E}(0)=(0,1/2] \ \ \text{and} \ \ \mathcal{E}(1)=[1/2,1).
$$
The following result is concerned about the almost surely finiteness of the random variable $\sup_{t\in \mathcal{E}(S)}\left(\chi_{\vk{b}}^2(t)- g(t)\right)$. 
\BT\label{Addition}
\aH{Let $\{X(t),t\in E\}$ be given as in  \netheo{Th0} with $E=(0,1)$. 
If $|f(S)|<\IF$ and conditions  \textbf{A}(S), \textbf{D}(S) are satisfied, then
 $$
 \pk{\sup_{t\in \mathcal{E}(S)}\left(\chi_{\vk{b}}^2(t)- g(t)\right)<\IF }= 1.
 $$
 If $|f(S)|=\IF$ and conditions  \textbf{A}(S), \textbf{B}(S) and \textbf{E}(S) are satisfied, then
 $$
 \pk{\sup_{t\in \mathcal{E}(S)}\left(\chi_{\vk{b}}^2(t)- g(t)\right)<\IF }= 1 \ \ \text{or} \ \ 0
 $$
as the integral $I_g(S)<\IF$ or $=\IF$.}
\ET
 The technical proof of \netheo{Addition} is presented in the Appendix.
 
  \aH{Next, we present our principle result. } % concerning the asymptotics of \eqref{chig}.}
\BT\label{TH1}
Let $\{X(t),t\in E\}$ be given as in  \netheo{Th0} with $E=(0,1)$. Then for each of the following scenarios we have that \eqref{PRC} holds for $E=(0,1)$.\\
  (i). $f(0)=-\IF, f(1)=\IF$, and conditions \textbf{A}(0), \textbf{B}(0), \textbf{C}(0), \textbf{A}(1), \textbf{B}(1), \textbf{C}(1) are satisfied;\\
  (ii). $f(0)=-\IF, f(1)<\IF$, and conditions \textbf{A}(0), \textbf{B}(0), \textbf{C}(0), \textbf{D}(1) are satisfied;\\
  (iii). $f(0)>-\IF, f(1)=\IF$, and conditions \textbf{D}(0),   \textbf{A}(1), \textbf{B}(1), \textbf{C}(1) are satisfied;\\
  (iv). $f(0)>-\IF, f(1)<\IF$, and conditions \textbf{D}(0), \textbf{D}(1) are satisfied.

\ET
\BRM\label{rem}  %a) In light of the proof of Theorem \ref{TH1}, we know that (\ref{K}) is needed only for $k=1$.\\
 Note that \netheo{TH1} can be easily extended to the case where $\cH{E}:=(0,T)$, with $T\in(0,\IF]$, by using a time-scaling. More precisely, suppose that the process $X$ and the trend function $g(\cdot)$ are well-defined on $E$. If $T$ is finite, we consider $Y(t)=X(Tt),t\in(0,1)$; if $T=\IF,$
we  consider $Z(t)=X(h(t)), t\in(0, 1)$,  %so that the problem over $(0,\IF)$ is transformed into a problem over $(0,1)$. %The corresponding conditions and results can be given with little effort so we leave them to readers.
where $h(t), t\in (0,1)$ is a monotone function with $\lim_{t\rw 0}h(t)=\IF$ and $\lim_{t\rw 1}h(t)=0$. For example,
\BQNY
h(t)=\left\{\begin{array}{cc}
1/t,&0<t\leq 1/2,\\
4(1-t),&1/2\leq t<1.
\end{array}
\right.
\EQNY
Under analogue conditions as in \netheo{TH1} on the processes $Y$ or $Z$, we can obtain a similar result as in \eqref{PRC} for these two cases.
\ERM
%%%%%%%%%%%%%%%%%%%%%%%%%%%%%%%%%%%%5
\COM{For the discussion of the  asymptotic tail probability of $L_E^\nu$ defined in \eqref{INBR}. We present next a variant of the law of iterated logarithm for Brownian bridge. % is given in term of chi-square process with trend.

\BEL\label{lemBB} Let $\{B(t),t\in(0,1)\}$ be the standard Brownian bridge, and let $g(t), t\in(0,1)$ be a nonnegative continuous function satisfying $g(t)\uparrow \IF$ as $t\rw 0$ and $t\to1$. We have\\
i). If $\int_0^1\frac{\sqrt{g(t)}}{t(1-t)}e^{-\frac{g(t)}{2}}dt<\IF$, then
\BQNY
\pc{\sup_{t\in(0,1)}\left(\frac{B^2(t)}{t(1-t)}-g(t)\right)<\IF}=1;
\EQNY
ii). If $\int_0^1\frac{\sqrt{g(t)}}{t(1-t)}e^{-\frac{g(t)}{2}}dt=\IF$, then
\BQNY
\pc{\sup_{t\in(0,1)}\left(\frac{B^2(t)}{t(1-t)}-g(t)\right)<\IF}=0.
\EQNY
\EEL
The proof of the above lemma follows directly from \nelem{A2} by using the fact that $B(t)=(1-t)W(\frac{t}{1-t}),t\in(0,1)$, with $W$ the standard Brownian motion.}
%%%%%%%%%%%%%%%%%%%%%%%%%%%%%%%%%%%%%%%%%%%%%%%%%%%%%%%5

\COM{When $\abs{f(S)}=\IF$, it is of interest to see the sufficiency and necessity of condition {\bf C}. For this purpose we further impose
%Additionally to conditions {\bf A-C}, we impose
 the following condition:}

As an application of \aH{\netheo{Addition} and Theorem \ref{TH1}}, we obtain the following result concerning
the supremum of the squared chi-square normalized standard Brownian bridge with trend, $L_{E}^\nu$ defined in \eqref{INBR}.
\BK\label{cor}
Let  $L_{E}^\nu$ be defined in \eqref{INBR} with $E=(0,1)$. We have:

 If $\nu>3/4$, then,   as $u\to\IF$
\BQNY
\pk{L_{E}^\nu >u}\sim\frac{\sqrt{u}e^{-u/2}}{\sqrt{2\pi}}\int_{0}^1\frac{1}{t(1-t)}e^{-g_\nu(t)}dt;
\EQNY
If $\nu\leq 3/4$, then
\BQNY
\pk{L_{E}^\nu=\IF}=1.
\EQNY
\EK

\COM{\BRM\label{rem:nu}
We see from the proof of the last theorem that $\int_0^{1}\frac{(g_\nu(t))^{1/2}}{t(1-t)}e^{-g_\nu(t)}dt<\IF$  if and only if $\nu>3/4$. Similarly, we can show that
$\int_{0}^1\frac{1}{t(1-t)}e^{-g_\nu(t)}dt<\IF$ if and only if $\nu>1/2$. Clearly, for any $\nu\in(1/2,3/4]$, the former integral is infinite but the latter is finite. This observation indicates that our condition {\bf C}  in \netheo{TH1} may be optimal; it cannot be weakened to $J_E^g<\IF$ as discussed in Section 2.
\ERM}

%%%%%%%%%%%%%%%%%%%%%%%%%%%%%%%%%%%%%%%%%%%%%%%%555
\COM{
Comparing the integral $J_g$ defined in (\ref{Jg}) appearing in Theorem \ref{TH1} to the condition {\bf C(S)} required in Theorem \ref{TH1}, we find that it seems the condition {\bf C(S)} requires more. The natural question is that can we relax {\bf C(S)} to $J_g<\IF$ ? More precisely, if we replace {{\bf C(S)}} with
$\abs{\int_{T/2}^S (C(t))^{1/\alpha}e^{-\frac{g(t)}{2}}dt}<\IF$, does Theorem \ref{TH1} still hold? Corollary \ref{cor2} suggests that for $n=1$ and $X$ to be the standard Brownian motion,  {\bf C(S)} is the optimal condition to derive the asymptotics. For example, if let $$g(t)=2\ln\ln(e^2/t)+2\rho\ln\ln\ln(e^3/t),$$
then $$\frac{1}{t}(g(t))^{\frac{1}{2}}e^{-\frac{g(t)}{2}}\sim \frac{ \mathcal{C }  }{t\ln(e^2/t)(\ln\ln(e^3/t))^{\rho-1/2}},$$
and $$\frac{1}{t}e^{-\frac{g(t)}{2}}\sim \frac{ \mathcal{C }  _1}{t\ln(e^2/t)(\ln\ln(e^3/t))^{\rho}}.$$

Thus we know that if $\rho\in(1,3/2]$, $J_g<\IF$ for $E=(0,1)$, whereas  $\int_{0}^1 \frac{1}{t}(g(t))^{\frac{1}{2}}e^{-\frac{g(t)}{2}}dt=\IF$. For this case, (\ref{e1}) doesn't hold since $\pc{\sup_{t\in(0,1]}\frac{B^2(t)}{t}-g(t)=\IF}=1$.
}%%%%%%%%%%%%%%%%%%%%%%%%%%%%%%%%%%%%%%%%55555

%Let $\{B_H(t), t\ge0\}$  be the standard fractional Brownian motion (fBm) with Hurst index $H\in(0,1)$ and covariance function given by
%$$
%\Cov(B_H(s),B_H(t))=\frac{1}{2}(t^{2H}+s^{2H}-\abs{t-s}^{2H}),\ \ \ s,t\ge0.
%$$
Given  wide applications of fBm in various fields, we give below a result concerning the tail asymptotics of the squared  normalized fBm with trend.  % one type of iterated logarithm law for Brownian motion.

 \BK\label{cor2} Let $\{X_i(t),t>0\}, i\leq n$ in (\ref{INCHI}) be independent copies of $\{B_H(t)/t^H, t>0\}$, with %  which is  the standard fBm with Hurst index
 $H\in(0,1)$, and let $g(\cdot)$ be a nonnegative continuous function on (0,1] such that $g(t)\uparrow \IF$ as $t\rw 0$. We have: \\
If $\int_{0}^1 \frac{1}{t}(g(t))^{k/2+1/(2H)- 1}e^{-\frac{g(t)}{2}}dt<\IF$, then,   as $u\to\IF$
\BQNY
\pc{\sup_{t\in(0,1]}\LT(\frac{\chi_{\vk{b}}^2(t)}{t^{2H}}-g(t)\RT)>u}\sim \mathcal{H}_{2H}\prod_{i=k+1}^n\left(1-b_i^2\right)^{-1/2}\frac{u^{\frac{k}{2}+\frac{1}{2H}-1}e^{-u/2}}{2^{\frac{k}{2}+\frac{1}{2H}-1}\Gamma(k/2)}\int_0^1\frac{1}{t}e^{-\frac{g(t)}{2}}dt; %\ \ u\rw\IF.
\EQNY
If $\int_{0}^1 \frac{1}{t}(g(t))^{k/2+1/(2H)- 1}e^{-\frac{g(t)}{2}}dt=\IF$, then
\BQNY
\pc{\sup_{t\in(0,1]}\LT(\frac{\chi_{\vk{b}}^2(t)}{t^{2H}}-g(t)\RT)=\IF}=1.
\EQNY
\EK

 Let $\vk{W}_n(t)=(W_1(t), W_2(t),\cdots, W_n(t)), t\ge0, n\geq 1$ be the standard $n$-dimensional Brownian motion, and let $||\vk{W}_n(t)||:=\sqrt{\sum_{i=1}^n{W_i^2(t)}}, t\geq 0$ be the Bessel process of order $n$. As a special case of Corollary \ref{cor2}, \aH{we have the following result. }% which  is closely related to the Kolmogorov-Dvoretsky-Erd\H{o}s integral test of the Bessel process; see, e.g.,  \aH{page 163 in \cite{IM1974} and Theorem A in} \cite{KLS96} and the references therein.
 \BK\label{cor4}
Let $g(\cdot)$ be a nonnegative continuous function on (0,1] such that $g(t)\uparrow \IF$ as $t\rw 0$. We have:

 If $\int_0^1\frac{1}{t}(g(t))^{n/2}e^{-\frac{g(t)}{2}}dt<\IF$, then,   as $u\to\IF$
  \BQN\label{CE1}
\pk{ \sup_{t\in (0,1]}\left(\frac{||\vk{W}_n(t)||^2}{t}- g(t)\right)>u}\sim \frac{2^{-n/2}u^{n/2}e^{-u/2}}{\Gamma(n/2)}\int_0^1\frac{1}{t}e^{-\frac{g(t)}{2}}dt;
\EQN
 If $\int_0^1\frac{1}{t}(g(t))^{n/2}e^{-\frac{g(t)}{2}}dt=\IF$, then
  \BQN\label{CE2}
\pk{ \sup_{t\in (0,1]}\left(\frac{||\vk{W}_n(t)||^2}{t}- g(t)\right)=\IF}=1.
\EQN
 \EK
 \aH{\BRM
 It can be shown that the classical law of interacted logarithm for Bessel process follows from   Corollary 2.8. In fact, 
define
 $$
 g_{\rho}(t)=2\ln\ln(e^2/t)+2\rho\ln\ln\ln(e^3/t),\ t\in(0,1].
 $$
 It follows that
 \BQNY
 \frac{1}{t}(g_\rho(t))^{n/2}e^{-\frac{g_\rho(t)}{2}}\sim \frac{ \mathcal{C }  }{t\ln(1/t)(\ln\ln(1/t))^{\rho-n/2}} %,\ \ \frac{1}{t}e^{-\frac{g_\rho(t)}{2}}\sim \frac{ \mathcal{C }  _1}{t\ln(1/t)(\ln\ln(1/t))^{\rho}}
 \EQNY
as $t\rw 0$, with $ \mathcal{C }  $  some positive constant. Elementary calculations show that  $\int_0^1\frac{(g_\rho(t))^{n/2}}{t}e^{-\frac{g_\rho(t)}{2}}dt<\IF$  holds if and only if $\rho>1+n/2$. Then,  by Corollary \ref{cor4} one can show that, for any $\rho>1+n/2$
\BQNY
\limsup_{t\to 0} \frac{||\vk{W}_n(t)||^2}{t g_\rho (t)}\le 1,
\EQNY
and
 for any $0< \rho\le 1+n/2$
\BQNY
\limsup_{t\to 0} \frac{||\vk{W}_n(t)||^2}{t g_\rho (t)}\ge 1.
\EQNY
Consequently, we arrive at the classical law of interacted logarithm as follows:
\BQNY
\limsup_{t\to 0} \frac{||\vk{W}_n(t)||^2}{2t \ln\ln(1/t) }=1.
\EQNY
\ERM}
\COM{  \BRM\label{rem:rho}
Define the following function
 $$g_{\rho}(t)=2\ln\ln(e^2/t)+2\rho\ln\ln\ln(e^3/t),\ t\in(0,1].$$
 It follows that
 \BQNY
 \frac{1}{t}(g_\rho(t))^{n/2}e^{-\frac{g_\rho(t)}{2}}\sim \frac{ \mathcal{C }  }{t\ln(1/t)(\ln\ln(1/t))^{\rho-n/2}} %,\ \ \frac{1}{t}e^{-\frac{g_\rho(t)}{2}}\sim \frac{ \mathcal{C }  _1}{t\ln(1/t)(\ln\ln(1/t))^{\rho}}
 \EQNY
as $t\rw 0$, with $ \mathcal{C }  $  a positive constant. Elementary calculations show that condition {\bf C}$(0)$ (which in our case is equivalent to $\int_0^1\frac{(g_\rho(t))^{n/2}}{t}e^{-\frac{g_\rho(t)}{2}}dt<\IF$) holds if and only if $\rho>1+n/2$.
 %and that $J_{(0,1)}^{g_\rho}=\int_{0}^1 \frac{1}{t}e^{-\frac{g_\rho(t)}{2}}dt<\IF$ holds if and only if $\rho>1$, which implies that for any $\rho\in(1,1+n/2]$, we have $J_{(0,1)}^{g_\rho}<\IF$  but (\ref{CE1}) does not hold. This observation indicates that our condition {\bf C}  in \netheo{TH1} may be optimal; it cannot be weakened to $J_{(0,1)}^g<\IF$ as discussed above.
 \ERM
}

\section{Further Results and Discussions}
In this section, we discuss an extension of (\ref{PRC})  to the case where $X_i$'s are not identically distributed. Namely, we assume that $X_i,1\leq i\leq n$ are independent locally stationary Gaussian processes on $(0,1)$  with a.s. continuous sample paths, unit variance and correlation functions $r_i(\cdot,\cdot)$ satisfying for any $1\leq i\leq k(\le n),$
\BQN\label{cos}
\lim_{h\to 0}\frac{1-r_i(t,t+h)}{K^2(\abs{h})}=C_i(t),
\EQN
and, for any $k+1\leq i\leq n$,
\BQN\label{cos1}
\lim_{h\to 0}\frac{1-r_i(t,t+h)}{K_i^2(\abs{h})}=C_i(t)
\EQN
hold uniformly in $t\in I$, any compact interval in $(0,1)$,
where $K(\cdot)$ and $K_i(\cdot), k+1\leq i\leq n$ are regularly varying functions at $0$ with index $0<\alpha/2< \alpha_{k+1}/2\leq \cdots \leq \alpha_n/2< 1$ respectively, and $C_i(\cdot), 1\leq i\leq n$ are positive continuous functions over $(0,1)$ satisfying   $C_i(0)=\IF$ or $C_i(1)=\IF$.

For simplicity we shall assume $\vk{b}=\vk{1}=(1,1,\cdots,1)\in\R^n$, and consider the asymptotics of
$$
\pk{\sup_{t\in E}\Bigl(\chi_{\vk{1}}^2(t)-g(t)\Bigr)>u}
$$
as $u\to\IF$, for certain admissible functions $g(\cdot)$ and $E$ either $(0,1)$ or a compact interval in $(0,1)$.

First, we present the result with $E$ being a  compact interval in $(0,1)$. Recall that $q(u)=\overleftarrow{K}(u^{-1/2})$.

\BT\label{Thd0} Let $X_i,1\leq i\leq n$ be defined as above with
correlation functions satisfying \eqref{cos} and \eqref{cos1}. If further $\max\{r_i(s,t), 1\leq i\leq n\}<1$ holds for all $s\neq t \in E$, a compact interval of $(0,1)$, then for any nonnegative continuous trend function $g(\cdot)$
we have, as $u\to \IF$
\BQN\label{dmr}
&&\pk{\sup_{t\in E}\Bigl(\chi_{\vk{1}}^2(t)-g(t)\Bigr)>u}\nonumber\\
&&\sim (2\pi)^{-n/2}\mathcal{H}_{\alpha}
\eeH{\frac{u^{n/2-1}e^{-u/2}}{q(u)}}\\
&&\times\int_{(t,\vk{\theta}) \in D_E }\left(C_1(t)\prod_{j=2}^n\cos^2(\theta_j)+\sum_{i=2}^kC_i(t)\sin^2(\theta_i)
\prod_{j={i+1}}^n\cos^2(\theta_j)\right)^{1/\alpha}\prod_{i=3}^n(\cos(\theta_i))^{i-2}
e^{-\frac{g(t)}{2}}dtd\vk{\theta},\nonumber
\EQN
\eeH{where $\vk{\theta}=(\theta_2,\cdots,\theta_n)$ and $D_E=E\times[-\pi,\pi]\times[-\pi/2,\pi/2]^{n-2}$}.
\ET
Next, we consider the case $E=(0,1)$. Similarly as in Section 2, we introduce a crucial function
\BQNY
f^*(t)=\int_{1/2}^tC^*(s)ds,\ \ \ t\in(0,1),
\EQNY
with $C^*(t)=\max\{C_i^{1/\alpha}(t), 1\leq i\leq k, C_i^{1/{\alpha_i}}(t), k+1\leq i\leq n\}, t\in(0,1)$.
We denote by $\invfs(t), t\in(f^*(0), f^*(1))$ the inverse function of $f^*(t),t\in(0,1)$. Further, for any $d>0$,
 let $s_{j,d}^{(1)}=\invfs(jd)$, $j\in \mathbb{N}\cup\{0\}$ if $f^*(1)=\IF$, and let $s_{j,d}^{(0)}=\invfs(-jd)$, $j\in \mathbb{N}\cup\{0\}$ if $f^*(0)=-\IF$. Denote $\Delta_{j,d}^{(1)}=[s_{j-1,d}^{(1)},s_{j,d}^{(1)}], j\in \mathbb{N}$ and
  $\Delta_{j,d}^{(0)}=[s_{j,d}^{(0)},s_{j-1,d}^{(0)}], j\in \mathbb{N}$, which give a partition of $[1/2,1)$ in the case $f^*(1)=\IF$ and a partition of $(0,1/2]$ in the case $f^*(0)=-\IF$, respectively.
 % then we can derive a division of $[T/2,T)$ denoted by $\Delta_{j,d}^{(T)}=[s_{j-1,d}^{(T)},s_{j,d}^{(T)}], j\in \mathbb{N}$. Similarly, if, we can derive a division for, i.e., $\Delta_{j,d}^{(0)}=[s_{j,d}^{(0)},s_{j-1,d}^{(0)}], j\in \mathbb{N}$ with.
 % We first consider the case that
%$$\cH{I(E):=\int_{t\in E}}(C(t))^{1/\alpha}dt=\IF.$$

  Let $S\in\{0,1\}$. Additionally to condition {\bf A}$(S)$ and analogously to conditions {\bf B}$(S)$, {\bf C}$(S)$, {\bf D}$(S)$ we impose the following (scenario-dependent) restrictions on the trend function $g(\cdot)$ and the correlation functions $r_i(\cdot,\cdot)$'s.\\
%\EH{$\lim_{t\rw \cH{S} }g(t)=\IF$} and there exist a constant $0<\delta_g<T/2$  such that $g$ is monotone over $E\cap (S-\delta_g, S+\delta_g)$ \cH{for $S\in \{0,T\}$}. \\
\textbf{Condition B'}$(S)$: Suppose that there exists some constant   $d_0>0$ such that
\BQNY
\limsup_{j\rw \IF}\sup_{t\neq s \in\Delta^{(S)}_{j,d_0}}\frac{1-r_i(t,s)}{K^2(|f^*(t)-f^*(s)|)}\cH{ < \IF},\ \ \forall 1\le i\le n.
 \EQNY
\textbf{Condition C'}$(S)$:  It holds  that
 \BQNY
  \abs{\int_{1/2}^S C^*(t)\frac{(g(t))^{\frac{n}{2}-1  }}{q(g(t))}e^{-\frac{g(t)}{2}}dt}<\IF.
 %~~ \int_{0}^{d\delta_0}\widetilde{R}_{d}(t)(1+g(t))^{(\frac{k}{2}-1)_{+} +1+\eta_1}e^{-\frac{g(T-t)}{2}}dt<\IF.
 \EQNY
 \textbf{Condition D'}($S$): The following is satisfied  % For any $S\in\{0,T\}$   %There exists a constant $M_F$ such that
\BQNY
\limsup_{\delta\rw 0}\sup_{t\neq s\in (0,\delta)}\frac{1-r_i(|S-t|, |S-s|)}{K^2(|f^*(|S-t|)-f^*(|S-s|)|)} \cH{< \IF},\ \ \forall 1\le i\le n.
\EQNY

We present below the main result of this section.

\BT\label{TH2}
Let $X_i, 1\leq i\leq n$ be given as in  \netheo{Thd0}. Then, for each of the following scenarios we have that \eqref{dmr} holds for $E=(0,1)$.\\
(i). $f^*(0)=-\IF, f^*(1)=\IF$, and conditions \textbf{A}(0), \textbf{B'}(0), \textbf{C'}(0), \textbf{A}(1), \textbf{B'}(1), \textbf{C'}(1) are satisfied;\\
(ii). $f^*(0)=-\IF, f^*(1)<\IF$, and conditions \textbf{A}(0), \textbf{B'}(0), \textbf{C'}(0), \textbf{D'}(1) are satisfied;\\
(iii). $f^*(0)>-\IF, f^*(1)=\IF$, and conditions \textbf{D'}(0),   \textbf{A}(1), \textbf{B'}(1), \textbf{C'}(1) are satisfied;\\
(iv). $f^*(0)>-\IF, f^*(1)<\IF$, and conditions \textbf{D'}(0), \textbf{D'}(1) are satisfied.

\ET
\BRM
%a) If $K^2_i(t)=d_i t^{\alpha}(1+o(1)), t\to 0$ for some positive constants $d_i$, then $\eta$ in condition {\bf C'} can be taken as $ 0$.\\
%ii) We note that if $n=k$ and $X_i, 1\leq i\leq k$ are iid, then Theorem \ref{Thd0} and Theorem \ref{TH2} are reduced to Theorem \ref{Th0} and Theorem \ref{TH1} with $n=k$.\\
%iii) Actually, we can also choose $C^*(t)=\sum_{i=1}^k C_i^{1/\alpha}(t)+\sum_{i=k+1}^nC_i^{1/{\alpha_i}}$ but then we need larger $g$.\\
 The problem becomes more difficult  when $\vk{b}\neq \vk{1}$; the difficulty comes from the fact that the expansion of the  correlation function of $Y_{\vk{b}}(t,\vk{\theta})$ as in (\ref{PRC2}) is too complicated.
\ERM

As a direct application of \netheo{TH2}, we obtain the following result concerning the tail asymptotics of the supremum of a chi-square process with trend.

\BK\label{cor3}
 Let $\{B(t),t\in[0,1]\}$ be the  standard Brownian bridge, $\{W(t), t\in [0,1]\}$ be the standard Brownian motion, and $\{B_H(t),t\in[0,1]\}$ be the standard fBm with Hurst index $H\in(1/2,1)$. Further, let $g(\cdot)$ be a nonnegative continuous function on $(0,1)$. If $g(\cdot)$ satisfies {\bf A}$(S)$ with $S\in\{0,1\}$  and $\int_0^1\frac{1}{t(1-t)}(g(t))^{3/2}e^{-\frac{g(t)}{2}}dt<\IF$, then, as $u\to\IF$
\BQNY
\pc{\sup_{t\in(0,1)}\left(\frac{B^2(t)}{t(1-t)}+\frac{W^2(t)}{t}+\frac{B_H^2(t)}{t^{2H}}-g(t)\right)>u}\sim \frac{u^{3/2}e^{-u/2}}{3\sqrt{2\pi}}\int_0^1\frac{2-t}{t(1-t)}e^{-\frac{g(t)}{2}}dt.
\EQNY
\EK

\section{Proofs}
For notational simplicity we denote $\varphi(u)=\prod_{i=k+1}^n\left(1-b_i^2\right)^{-1/2}(q(u))^{-1}u^{k/2-1}e^{-u/2}$ (recall $q(u)=\overleftarrow{K}(u^{-1/2})$).

In what follows, we use $ \mathcal{C }  , \mathcal{C }  _1, \mathcal{C }  _2,\cdots$ to denote unspecified positive and finite constants which may not be the same from line to line. % More specific constants are numbered as $$.

\prooftheo{Th0} Using the classical approach when dealing with extremes of chi-square processes we reduce  the problem to the study of extremes of Gaussian random fields; see e.g., \cite{JAPIT, Pit96}. To this end, we introduce two particular Gaussian random fields, namely  (denote $D=D_E=E\times[-\pi,\pi]\times[-\pi/2,\pi/2]^{n-2})$
\BQN
Y_{\vk{b}}(t,\vk{\theta})&=&\sum_{i=1}^nb_iX_i(t)v_i(\vk{\theta}),\label{Y}\\
Y_{\vk{b}}^{(u)}(t,\vk{\theta})&=&\sum_{i=1}^nb_iX_i^{(u)}(t)v_i(\vk{\theta}), \ \ \ (t,\vk{\theta})=(t,\theta_2,\cdots,\theta_n)\in D,\nonumber
\EQN
where $X_i^{(u)}(t)=\frac{X_i(t)}{\sqrt{1+g(t)/u}}$, $t\in E$, and $v_n(\vk{\theta})=\sin( \theta_n), v_{n-1}(\vk{\theta})=\sin (\theta_{n-1})\cos(\theta_n),\cdots, v_1(\vk{\theta})=\cos (\theta_n)\cdots\cos( \theta_2)$ are spherical coordinates. In view of  \cite{Pit96},  for any $u>0$
\BQNY
\mathbb{P}\left(\sup_{t\in E}\eeH{\bigl (}\chi_{\vk{b}}^2(t)-g(t) \eeH{\bigr )}>u\right)=\mathbb{P}\left(\sup_{(t,\vk{\theta})\in \eeH{D}}Y_{\vk{b}}^{(u)}(t,\vk{\theta})>\sqrt{u}\right).
\EQNY
 Let $A=E\times[-\pi,\pi]\times[-\pi/2,\pi/2]^{k-2}$ and $B_u=[-m(u),m(u)]^{n-k}$ with $m(u)= {\ln (u)}/{\sqrt{u}}$.
 \eeH{For any $u>0$}
 \BQNY
 \pi(u)\leq \mathbb{P}\left(\sup_{(t,\vk{\theta})\in D}Y_{\vk{b}}^{(u)}(t,\vk{\theta})>\sqrt{u}\right)\leq \pi(u)+\mathbb{P}\left(\sup_{(t,\vk{\theta})\in D\setminus (A\times B_u)}Y_{\vk{b}}^{(u)}(t,\vk{\theta})>\sqrt{u}\right),
 \EQNY
 where
 $$\pi(u):=\mathbb{P}\left(\sup_{(t,\vk{\theta})\in A\times B_u}Y_{\vk{b}}^{(u)}(t,\vk{\theta})>\sqrt{u}\right).$$

\EH{Since} the variance function of $Y_{\vk{b}}^{(u)}$ satisfies \eeH{for $u>0$}
\BQNY
%\left(\E{\left(Y_{\vk{b}}^{(u)}(t,\vk{\theta})\right)^2}\right)^{1/2}&=&\sqrt{\frac{u}{u+g(t)}}\left(1-(1-b_n^2)\sin^2\theta_n-%\sum_{i=k+1}^{n-1}
%(1-b_i^2)\left(\prod_{j=i+1}^{n}\cos^2\theta_j\right)\sin^2 \theta_i\right)^{1/2},
\E{\left(Y_{\vk{b}}^{(u)}(t,\vk{\theta})\right)^2} = \frac{u}{u+g(t)}\left(1-(1-b_n^2)\sin^2(\theta_n)-\sum_{i=k+1}^{n-1}
(1-b_i^2)\left(\prod_{j=i+1}^{n}\cos^2(\theta_j)\right)\sin^2 (\theta_i)\right),
\EQNY
 %as $\theta_i\rw 0, k+1\leq i\leq n.$ \\
 %\eeH{Question: Why $\theta \to 0$? Above there is equality, there is no $(1+ o(1))$ term.} This should be equality for all $\theta$!!!\\
 \eeH{for all $u$ large} \EH{we have}
  %\eeH{Question: What does $D-A\times B_u$ means, is it $D\setminus (A\times B_u)$?}\\
  %\eeH{Question: The constant $ \mathcal{C }  $ two lines above should be different from that below! Some remarks is needed.}\\
  \BQNY
  \sup_{(t,\vk{\theta})\in D\setminus (A\times B_u)}\LT(\E{(Y_{\vk{b}}(t,\vk{\theta}))^2}\RT)^{1/2}\leq \left(1-(1-b_n^2)\sin^2 (m(u))\right)^{1/2}\leq 1- \mathcal{C }  m^2(u).
  \EQNY
%with some constant $ \mathcal{C }  >0$.
Further, it can be shown from (\ref{INCO}) that %Since  Due to
 \BQNY
 \E{\left(Y_{\vk{b}}(t,\vk{\theta})-Y_{\vk{b}}(s,\vk{\theta}^{'})\right)^2}\leq  \mathcal{C }  _1|(t,\vk{\theta})-(s,\vk{\theta}^{'})|^{\alpha/2},~~(t,\vk{\theta}), (s,\vk{\theta}^{'})\in D.
 \EQNY
%with some constant $ \mathcal{C }  _1>0$.
Consequently, the Piterbarg \eeH{inequality} (see e.g., Theorem 8.1 in \eeH{\cite{Pit2001} or \cite{Pit96}}) implies
  \BQNY
  \mathbb{P}\left(\sup_{(t,\vk{\theta})\in D\setminus (A\times B_u)}Y_{\vk{b}}^{(u)}(t,\vk{\theta})>\sqrt{u}\right)&\leq& \mathbb{P}\left(\sup_{(t,\vk{\theta})\in D\setminus (A\times B_u)}Y_{\vk{b}}(t,\vk{\theta})>\sqrt{u}\right)\\
  &\leq& \mathcal{C }  _2u^{\frac{2n}{\alpha}}\Psi(\sqrt{u}/(1- \mathcal{C }  _3m^2(u)))\\
  &=&o(\pi(u)), ~~u\rw\IF,
  \EQNY
  %with $ \mathcal{C }  _2,  \mathcal{C }  _3$ two positive constants,
  \eeH{ where the last equality} is based on the following lower bound
  \BQNY
  \pi(u)\geq \mathbb{P}\left(Y_{\vk{b}}^{(u)}(\delta,0,\cdots,0)>\sqrt{u}\right)=\Psi(\sqrt{u+g(\delta)})(1+o(1)),
  ~~u\rw\IF.
  \EQNY
\eeH{Consequently,}
  \BQNY
  \mathbb{P}\left(\sup_{(t,\vk{\theta})\in D}Y_{\vk{b}}^{(u)}(t,\vk{\theta})>\sqrt{u}\right)=\pi(u)(1+o(1)), \ \ \ u\rw\IF.
  \EQNY
To complete the proof, it is thus sufficient to focus on the asymptotics of $\pi(u)$ as $u\to\IF.$ \eeH{Next, we  split} the rectangles $A$ and $\widetilde{A}=E \times[-\pi+\delta_1,\pi-\delta_1]\times[-\pi/2+\delta_1,\pi/2-\delta_1]^{k-2}$ with $\delta_1\in(0,\pi/2)$, into several subrectangles  denoted by $\{A_j\}_{j\in \Upsilon_1}$ and $\{\widetilde{A}_j\}_{j\in\Upsilon_2}$ respectively. \eeH{Further,} let
  $L$ and $\widetilde{L}$ represent their maximum length of edges \EH{of these} subrectangles, respectively.
  It follows from  Bonferroni's inequality that
  \BQNY
  \sum_{j\in\Upsilon_2}\mathbb{P}\left(\sup_{(t,\vk{\theta})\in \widetilde{A}_j\times B_u}Y_{\vk{b}}^{(u)}(t,\vk{\theta})>\sqrt{u}\right)-\eeH{\Lambda}(u)\leq\pi(u)\leq\sum_{j\in\Upsilon_1}\mathbb{P}\left(\sup_{(t,\vk{\theta})\in A_j\times B_u}Y_{\vk{b}}^{(u)}(t,\vk{\theta})>\sqrt{u}\right),
  \EQNY
  where
  $$\eeH{\Lambda}(u)=\sum_{j<j_1\in\Upsilon_2}\mathbb{P}\left(\sup_{(t,\vk{\theta})\in \widetilde{A}_j\times B_u}Y_{\vk{b}}^{(u)}(t,\vk{\theta})>\sqrt{u},
  \sup_{(t,\vk{\theta})\in \widetilde{A}_{j_1}\times B_u}Y_{\vk{b}}^{(u)}(t,\vk{\theta})>\sqrt{u}\right).$$
  For any fixed $j$, we have
  \BQNY
  \mathbb{P}\left(\sup_{(t,\vk{\theta})\in A_j\times B_u}Y_{\vk{b}}^{(u)}(t,\vk{\theta})>\sqrt{u}\right)\leq \mathbb{P}\left(\sup_{(t,\vk{\theta})\in A_j\times B_u}Y_{\vk{b}}(t,\vk{\theta})>\sqrt{u+g_j}\right),
  \EQNY
  where
  $g_j:=\eeH{\min_{\eeH{t}\in h\circ A_j}}g(t)$ and $h$ is a projection function defined by
  \BQNY
  h(x_1,\cdots,x_k)=x_1.
  \EQNY
  Further,
  \BQNY
\left(\E{\left(Y_{\vk{b}}(t,\vk{\theta})\right)^2}\right)^{1/2}=1-\sum_{i=k+1}^{n}\frac{1-b_i^2}{2}\theta_i^2(1+o(1)), \ \ \ \theta_i\rw 0, \quad k+1\leq i\leq n
\EQNY
and
\BQN\label{PRC2}
\Corr\left(Y_{\vk{b}}(t,\vk{\theta}), Y_{\vk{b}}(s,\vk{\theta}^{'})\right)&=&1-C(t_0)K^2(|t-s|)(1+o(1))-\sum_{i=2}^{k-1}\frac{1}{2}\left(\prod_{l=i+1}^k\cos^2 \theta_l\right)(\theta_i-\theta_i^{'})^2(1+o(1))\nonumber\\
&&-\sum_{i=k}^n\frac{b_i^2}{2}(\theta_i-\theta_i^{'})^2(1+o(1))
\EQN
as $t,s\rw t_0, |\theta_i-\theta_i^{'}|\rw 0, 2\leq i\leq k, \theta_i\rw 0, k+1\leq i\leq n.$ Next we introduce some useful notation. Let $\vk{\theta}_i=(\theta_i,\cdots,\theta_k), 2\leq i\leq k$, and $h_i, 2\leq i\leq k$ denote projection functions defined by
  \BQNY
  h_i(x_1,\cdots,x_k)=(x_i,\cdots,x_k),\ \ \ 2\leq i\leq k.
  \EQNY
For any $0<\epsilon<1$, there exist $u_0>0$ and $L_0>0$ such that for all $u>u_0$ and $0<L<L_0$ \EH{(recall that
$L$ is the maximum length of the edges of subrectangles $\{A_j\}_{j\in \Upsilon_1}$)}
\BQNY
 \mathbb{P}\left(\sup_{(t,\vk{\theta})\in A_j\times B_u}Y_{\vk{b}}(t,\vk{\theta})>\sqrt{u+g_j}\right)\leq\mathbb{P}\left(\sup_{(t,\vk{\theta})\in A_j\times B_u}\widetilde{Z}_\epsilon(t,\vk{\theta})>\sqrt{u+g_j}\right)
\EQNY
holds for all $j\in\Upsilon_1$, where $\widetilde{Z}_\epsilon(t,\vk{\theta})=\frac{Z_\epsilon(t,\vk{\theta})}{1+(1-\epsilon)\sum_{i=k+1}^{n}\frac{1-b_i^2}{2}\theta_i^2}$ and $Z_{\epsilon}(t,\vk{\theta})$ is a stationary Gaussian process with unit variance and correlation function $r_{Z_\epsilon}(t,\vk{\theta})$ satisfying
\BQNY
r_{Z_\epsilon}(t,\vk{\theta})=1-(1+\epsilon)C_jK^2(|t|)(1+o(1))-\sum_{i=2}^{k}\frac{e_{i,j}+\epsilon}{2}\theta_i^2(1+o(1))-(1+\epsilon)
\sum_{i=k+1}^n\frac{b_i^2}{2}\theta_i^2(1+o(1)) %\ \ \ (t,\vk{\theta})\rw 0,
\EQNY
as $(t,\vk{\theta})\rw \vk{0}$, with $C_j=\max_{t\in h\circ A_j}C(t)$,  $e_{k,j}=1$ and $e_{i,j}=\max_{\vk{\theta}_{i+1}\in h_{i+1}\circ A_j}\prod_{l=i+1}^k\cos^2 (\theta_l), 2\leq i\leq k-1$. The existence of such correlation function $r_{Z_\epsilon}$ can be confirmed by the Assertion on page 265 in \cite{HP99}, see also the reference mentioned therein.
 Therefore, by using similar arguments as in Theorem 3.2 in \cite{HAJI2014} (see also Theorem 8.2 in \cite{Pit96})  we can show that
 \BQN\label{PRI}
\lim_{u\rw\IF}\frac{\mathbb{P}\left(\sup_{(t,\vk{\theta})\in A_j\times B_u}\widetilde{Z}_\epsilon(t,\vk{\theta})>\sqrt{u+g_j}\right)}{\varphi(u)} = a(\epsilon)(2\pi)^{-k/2}\mathcal{H}_{\alpha}C_j^{1/\alpha}e^{-\frac{g_j}{2}}\prod_{i=2}^k(e_{i,j})^{1/2}\EH{mes(A_j)},
 \EQN
 where $a(\epsilon)\rw 1$ as $\epsilon\rw 0.$
Consequently % Further,
 \BQN\label{PRU}
 \lim_{\epsilon\rw 0}\lim_{L\rw 0}\limsup_{u\rw\IF}\frac{\pi(u)}{\varphi(u)}&\leq&(2\pi)^{-k/2}\mathcal{H}_{\alpha}
 \int_{(t,\EH{\vk{\theta}_2})\in A}(C(t))^{1/\alpha}e^{-\frac{g(t)}{2}}\prod_{i=3}^{k}(\cos\theta_i)^{i-2}dtd\theta_2\cdots d\theta_k\nonumber\\
 &=&\frac{2^{1-k/2}}{\Gamma(k/2)}\mathcal{H}_{\alpha} J_E^g. %\int_{\delta}^{T-\delta}(C(t))^{1/\alpha}e^{-\frac{g(t)}{2}}dt.
 \EQN
 Similarly, \EH{with $\widetilde{L}$ the maximum length of the edges of the subrectangles $\{\widetilde{A_j}\}_{j\in \Upsilon_1}$ we have}
 \BQN\label{PRL}
 \lim_{\delta_1\rw 0}\lim_{\epsilon\rw 0}\lim_{\widetilde{L}\rw 0}\liminf_{u\rw\IF}\frac{\sum_{j\in\Upsilon_2}\mathbb{P}\left(\sup_{(t,\vk{\theta})\in \widetilde{A}_j\times B_u}Y_{\vk{b}}^{(u)}(t,\vk{\theta})>\sqrt{u}\right)}{\varphi(u)}\geq \frac{2^{1-k/2}}{\Gamma(k/2)}\mathcal{H}_{\alpha}J_E^g. %\int_{\delta}^{T-\delta}C(t)^{1/\alpha}e^{-\frac{g(t)}{2}}dt.
 \EQN
We consider next  the asymptotic of $\Lambda(u):=\Lambda_1(u)+\Lambda_2(u)$ (to be specified later). For $u>0$  we have
 \BQNY
 &&\Lambda_1(u):=\sum_{\widetilde{A}_j\bigcap\widetilde{A}_{j_1}=\emptyset}\mathbb{P}\left(\sup_{(t,\vk{\theta})\in \widetilde{A}_j\times B_u}Y_{\vk{b}}^{(u)}(t,\vk{\theta})>\sqrt{u},
  \sup_{(t,\vk{\theta})\in \widetilde{A}_{j_1}\times B_u}Y_{\vk{b}}^{(u)}(t,\vk{\theta})>\sqrt{u}\right)\\
  &&\leq\sum_{\widetilde{A}_j\bigcap\widetilde{A}_{j_1}=\emptyset}\mathbb{P}\left(\sup_{(t,s,\vk{\theta},\vk{\theta}^{'})\in \widetilde{A}_j\times \widetilde{A}_{j_1}\times B_u\times B_u}\eeH{\Bigl(}Y_{\vk{b}}(t,\vk{\theta})+Y_{\vk{b}}(s,\vk{\theta}^{'})\eeH{\Bigl)}>2\sqrt{u}\right).
 \EQNY
Further, there exists some $ \delta_0\in(0,4)$ such that,
for all $\widetilde{A}_j\bigcap\widetilde{A}_{j_1}= \emptyset$  %and further
 \BQN\label{PRV}
 \E{\left(Y_{\vk{b}} (t,\vk{\theta})+Y_{\vk{b}} (s,\vk{\theta}^{'})\right)^2}\leq 4-\delta_0,\ \ \  (t,\vk{\theta},s,\vk{\theta}^{'})\in \widetilde{A}_j\times B_u\times\widetilde{A}_{j_1}\times B_u
 \EQN
holds for all large $u.$ Consequently, in the light of Borell-TIS inequality (see e.g., \cite{AdlerTaylor})
 \BQN\label{PRL1}
 \Lambda_1(u)\leq  \mathcal{C }    e^{\frac{-(2\sqrt{u}-a)^2}{2(4-\delta_0)}}=o(\varphi(u)), ~~u\rw\IF,
 \EQN
 where %$ \mathcal{C }  _4$ is a positive constant and
$a= \E{\sup_{(t,\vk{\theta})\in D}Y_{\vk{b}}(t,\vk{\theta})}<\IF$.\\
% \eeH{Note: We need to be careful with constants $ \mathcal{C }  $ they appear all over...}\\
 Further, we have
 \BQNY
 \Lambda_2(u)&:=&\sum_{\widetilde{A}_j\bigcap\widetilde{A}_{j_1}\neq\emptyset}\mathbb{P}\left(\sup_{(t,\vk{\theta})\in \widetilde{A}_j\times B_u}Y_{\vk{b}}^{(u)}(t,\vk{\theta})>\sqrt{u},
  \sup_{(t,\vk{\theta})\in \widetilde{A}_{j_1}\times B_u}Y_{\vk{b}}^{(u)}(t,\vk{\theta})>\sqrt{u}\right)\\
  &\leq&\sum_{\widetilde{A}_j\bigcap\widetilde{A}_{j_1}\neq\emptyset}\left[\mathbb{P}\left(\sup_{(t,\vk{\theta})\in \widetilde{A}_j\times B_u}Y_{\vk{b}}^{(u)}(t,\vk{\theta})>\sqrt{u}\right)+
  \mathbb{P}\left(\sup_{(t,\vk{\theta})\in \widetilde{A}_{j_1}\times B_u}Y_{\vk{b}}^{(u)}(t,\vk{\theta})>\sqrt{u}\right)\right.\\
  &&\left.-\mathbb{P}\left(\sup_{(t,\vk{\theta})\in (\widetilde{A}_j\cup\widetilde{A}_{j_1})\times B_u}Y_{\vk{b}}^{(u)}(t,\vk{\theta})>\sqrt{u}\right)\right].
 \EQNY
\EH{Along the same lines of the proof above}
 \BQNY
 \limsup_{u\rw \IF}\frac{\Lambda_2(u)}{\varphi(u)}\leq 3^n (a(\epsilon,\widetilde{L})-b(\epsilon,\widetilde{L}))\frac{2^{1-k/2}}{\Gamma(k/2)}\mathcal{H}_{\alpha} J_E^g, %\int_{\delta}^{T-\delta}(C(t))^{1/\alpha}e^{-\frac{g(t)}{2}}dt,
 \EQNY
 where $a(\epsilon, \widetilde{L}), b(\epsilon, \widetilde{L})\rw 1$ as $\widetilde{L}\rw 0$ and $\epsilon\rw 0$, which implies that
 \BQN\label{PRL2}
 \lim_{\epsilon\rw 0}\lim_{\widetilde{L}\rw 0}\limsup_{u\rw\IF}\frac{\Lambda_2(u)}{\varphi(u)}=0.
 \EQN
Consequently, we see from (\ref{PRL}), (\ref{PRL1}) and (\ref{PRL2}) that
 \BQNY
 \lim_{\delta_1\rw 0}\lim_{\epsilon\rw 0}\lim_{\widetilde{L}\rw 0}\liminf_{u\rw\IF}\frac{\pi(u)}{\varphi(u)}\geq \frac{2^{1-k/2}}{\Gamma(k/2)}\mathcal{H}_{\alpha}J_E^g, %\int_{\delta}^{T-\delta}C^{1/\alpha}(t)e^{-\frac{g(t)}{2}}dt,
 \EQNY
 which together with (\ref{PRU})  establishes  the proof. \QED
 %\BRM The reason why we consider the Gaussian process over the set $\widetilde{A}$ is to avoid the case that the coefficients of correlation function of $Y_b(t,\EH{\vk{\theta}})$ in(\ref{PRV}) \eeH{are} equal to zero and to eliminate the case that the variance in (\ref{PRV}) is equal to 4.
 %For both these two cases, we cannot continue our proof in the usual way. Thus taking into account $\widetilde{A}$ is necessary.
 %\ERM

We present next two results which are crucial to the proof of \netheo{TH1}.

 \BEL\label{LEV}
 Suppose that $V_{\vk{b}}^2(t)=\sum_{i=1}^nb_i^2V_i^2(t), t\in\R,$ with  $1=b_1=\cdots=b_k> b_{k+1}\geq \cdots\geq b_n>0$, is a stationary chi-square process, where $V_i, 1\leq i\leq n,$ are independent copies  of a stationary Gaussian process $\{V(t), t\in \R\}$  with unit variance, and correlation function $r(\cdot)$ satisfying $r(t)<1$ for all $t\neq 0$ and $1-r(t)\sim K^2(\abs{t})$ as $t\rw 0$, where $K(\cdot)$ is given as in Introduction. Then for any positive constants $a,b\in \R$ satisfying $a<b$ we have
 \BQNY
 \mathbb{P}\left(\sup_{t\in[a,b]}V_{\vk{b}}^2(t)>u\right)\sim \mathcal{H}_{\alpha} (b-a) \frac{2^{1-k/2}}{\Gamma(k/2)}\varphi(u)
 \EQNY
as $u\to\IF.$
 \EEL
 \textbf{Proof.}  As in the proof of \netheo{Th0} we can choose to work with the supremum of a Gaussian random fields. The proof follows by similar arguments as those in \cite{Piterbarg94}, and  % By using the similar approach as in the proof of \netheo{Th0}, we can derive the above conclusion.
thus we omit the details here.\QED

 % We present a weak version of Slepian's  inequality for chi-square processes in the following lemma.
 \BEL\label{LES} ({\bf weak  Slepian's  inequality})
Let $V_i(t),W_i(t), t\in A= [c,d] \subset \R, 1\le i\le n$ be independent centered Gaussian processes with a.s. continuous sample paths. If for all $1\le i\le n$,
$$\E{V^2_i(t)}=\E{W_i^2(t)}, t\in A \ \ \text{ and}\ \ \E{V_i(s)V_i(t)}\geq \E{W_i(s)W_i(t)}, s,t\in A
$$
are satisfied, then, we have, for all $u>0$,
 \BQNY
\pk{\sup_{t\in A} \sum_{i=1}^n V_i^2(t)>u}\leq 2^n \pk{\sup_{t\in A} \sum_{i=1}^n W_i^2(t)>u}.
 \EQNY

 \EEL

{\bf Proof.}  If $n=1$, then a direct application of Slepian's inequality yields
$$ \pk{ \sup_{t\in A} V_1^2(t)> u} \le 2 \pk{ \sup_{t\in A} V_1(t)> \sqrt{u}} \le
2 \pk{ \sup_{t\in A} W_1(t)> \sqrt{u}} \le  2 \pk{ \sup_{t\in A} W_1^2(t)> u}.$$

%%%%%%%%%%%%%%%%%%%%%%%%%%55555555555555555
 Below we shall show the claim for $n\ge 2$. Let therefore
 $\widetilde{X}(t,\vk{v})=\sum_{i=1}^n V_{i}(t)v_i, (t,\vk{v})\in A\times S_{n-1}$ and $\widetilde{Y}(t,\vk{v})=\sum_{i=1}^nW_{i}(t)v_i, (t,\vk{v})\in A\times S_{n-1}$ be the associated Gaussian random fields. % and suppose for simplicity that $n\ge 2$.
  Let $K_j,1\leq j\leq 2^n$ denote all the quadrants of $\R^n$ and define $\Delta_j= S_{n-1} \cap K_j, 1\le j\le 2^n$.
%So $\Delta_j=\{\vk{v}\in S_{n-1}, \inf_{i\in O_j} v_i\geq 0 ~\text{and}~\sup_{i\in \bar O_j }v_i\leq 0\}, 1\leq j\leq 2^n$
%where $O_j$ is a subset of $\{ 1 \ldot n\}$ and $\bar O_j$ denotes its complement with respect to $\{1 \ldot n\}$.
It follows that for any $1\leq j\leq 2^n$
 \BQNY
 \E{\left(\widetilde{X}(t,\vk{v})\right)^2}= \E{\left(\widetilde{Y}(t,\vk{v})\right)^2}, ~~(t,\vk{v})\in A\times\Delta_j,
 \EQNY
 and
 \BQNY
 \mathbb{E}\left(\widetilde{X}(t,\vk{v})\widetilde{X}(s,\vk{v}^{'})\right)\geq  \mathbb{E}\left(\widetilde{Y}(t,\vk{v})\widetilde{Y}(s,\vk{v}^{'})\right) \ \ \ (t,\vk{v}), (s,\vk{v}^{'})\in A\times\Delta_j
 \EQNY
 which is due to the fact that the sign of $\vk{v}$ is the same as that of $\vk{v}^{'}$ in this case. Thus,
 by applying \eeH{Slepian's} inequality  we obtain, for all $u>0$,
 \BQNY
 \mathbb{P}\left(\sup_{(t,v)\in A\times \Delta_j}\widetilde{X}(t,\vk{v})>u\right)\leq \mathbb{P}\left(\sup_{(t,\vk{v})\in A\times \Delta_j}\widetilde{Y}(t,\vk{v})>u\right),
 \EQNY
hence the proof follows easily. \QED

\prooftheo{TH1} We present only the proof for the case $(ii)$; the same arguments apply to other cases. First note that for any small $\delta>0$
\BQN\label{MAIN}
I_{\delta,1}(u)&:=&\mathbb{P}\left(\sup_{t\in[\delta, 1-\delta]}\left(\chi_{\vk{b}}^2(t)-g(t)\right)>u\right)\nonumber\\
&&\leq \mathbb{P}\left(\sup_{t\in(0, 1)}\left(\chi_{\vk{b}}^2(t)-g(t)\right)>u\right)\\
&&\leq I_{\delta,1}(u)+I_{\delta,2}(u)+I_{\delta,3}(u),\nonumber
\EQN
where
$$I_{\delta,2}(u):= \mathbb{P}\left(\sup_{t\in(0,\delta]}\left(\chi_{\vk{b}}^2(t)-g(t)\right)>u\right),\ \ I_{\delta,3}(u)=: \mathbb{P}\left(\sup_{t\in[1-\delta,1)}\left(\chi_{\vk{b}}^2(t)-g(t)\right)>u\right).$$
Since   from \netheo{Th0}
\BQNY
\lim_{\delta\to0}\lim_{u\to\IF}\frac{I_{\delta,1}(u)}{\varphi(u)}=\mathcal{H}_{\alpha}\frac{2^{1-k/2}}{\Gamma(k/2)} \int_{0}^{1}(C(t))^{1/\alpha}e^{-\frac{g(t)}{2}}dt
\EQNY
 it is sufficient to show that
\BQN\label{eq:I23}
\limsup_{\delta\to0}\limsup_{u\to\IF}\frac{I_{\delta,2}(u)}{\varphi(u)}=\limsup_{\delta\to0}\limsup_{u\to\IF}\frac{I_{\delta,3}(u)}{\varphi(u)}=0.
\EQN
 We first consider $I_{\delta,2}(u)$. % since the asymptotic of $I_{\delta,3}(u)$ can be derived similarly.
Without loss of generality, we choose $d$ as a parameter taking values in $\{2^{-m}d_0, m\in \mathbb{N}^+\}$. It is straightforward that
\BQN\label{MAINI}
I_{\delta,2}(u)&\leq& \sum_{j=N_{d,\delta}}^{\IF}\mathbb{P}\left(\sup_{t\in \Delta_{j,d}^{(0)}}\left(\chi_{\vk{b}}^2(t)-g(t)\right)>u\right)\nonumber\\
&\leq&\sum_{j=N_{d,\delta}}^{\IF}\mathbb{P}\left(\sup_{t\in \Delta_{j,d}^{(0)}}\chi_{\vk{b}}^2(t)>u+g(s_{j-1,d}^{(0)})\right),
\EQN
where $N_{d,\delta}=[-\frac{f(\delta)}{d}]$ with $[\cdot]$ the ceiling function.
Further, it follows from condition \textbf{B}(0)  that there exist some  $M_1>1, \delta_1>0$ and $d_1>0$ such that for all $\delta\in(0,\delta_1)$, $d\in(0,d_1)$, $j\geq N_{d,\delta}$,
\BQNY
1-r(t,s)\leq M_1 K^2(|f(t)-f(s)|)\leq \frac{1}{2}K^2((4M_1)^{1/\alpha}|f(t)-f(s)|), ~~t,s\in \Delta_{j,d}^{(0)}
\EQNY
hold. Next let $r(\cdot)$ be as  in \nelem{LEV}. Then for some constant $ \delta_2>0$
 \BQN\label{eqr}
  1-r(t)\geq \frac{1}{2}K^2(\abs{t}),\ \ \ t\in(0,\delta_2).
 \EQN
Therefore, for any   $d\in(0,\min(d_1, (4M_1)^{-1/\alpha}\delta_2)$ and $\delta\in(0,\delta_1)$
\BQNY
%\E{\left(X(t)\right)^2}&=&\E{\left(V((4M_1)^{1/\alpha}f(t))\right)^2},\\
\mathbb{E}\left(X(t)X(s)\right) \geq  \mathbb{E}\left(V((4M_1)^{1/\alpha}f(t))V(((4M_1)^{1/\alpha}f(s))\right)
\EQNY
holds for all $t, s\in\Delta_{j,d}^{(0)}$ and $j\geq N_{d,\delta}$, where $V$ is the stationary Gaussian process given in \nelem{LEV}.
%Apparently, we can derive the desired correlation inequality by the above inequality and (\ref{eqr}).
Consequently, in the light of \nelem{LES}  we have for $ \delta\in(0,\delta_1)$,  $d\in(0,\min(d_1, (4M_1)^{-1/\alpha}\delta_2)$  and $j\geq N_{d,\delta}$
\BQNY
\mathbb{P}\left(\sup_{t\in\Delta_{j,d}^{(0)}}\chi_{\vk{b}}^{2}(t)>u+g(s_{j-1,d}^{(0)})\right)&\leq& 2^n \mathbb{P}\left(\sup_{t\in \Delta_{j,d}^{(0)}}V_{\vk{b}}^{2}((4 M_1)^{1/\alpha}f(t))>u+g(s_{j-1,d}^{(0)})\right)\\
&\leq& 2^n\mathbb{P}\left(\sup_{t\in \Delta_d}V_{\vk{b}}^{2}(t)>u+g(s_{j-1,d}^{(0)})\right)
\EQNY
holds, where $\Delta_d=[0,(4M_1)^{1/\alpha}d]$ and $V_{\vk{b}}^{2}$  is the chi-square process given in \nelem{LEV}.
%Further
%\BQNY
%I_{\delta,2}(u)\leq 2^n\sum_{j=N_{d,\delta}}^{\IF}\mathbb{P}\left(\sup_{t\in \Delta_{d}}V^2_{\vk{b}}(t)>u+g(s_{j,d}^{(0)})\right).
%\EQNY
%With the aid of
Next, by \nelem{LEV}  we derive for any  $d\in(0,\min(d_1, (4M_1)^{-1/\alpha}\delta_2)$ and $\delta\in(0,\delta_1)$ %and $u$ is large enough, then
\BQN\label{eqs}
 &&\mathbb{P}\left(\sup_{t\in \Delta_d}V^2_{\vk{b}}(t)>u+g(s_{j-1,d}^{(0)})\right)\nonumber\\
&& \ \ \sim  H_\alpha\frac{2^{1-k/2}}{\Gamma(k/2)}(4M_1)^{1/\alpha}d\ \prod_{i=k+1}^n\left(1-b_i^2\right)^{-1/2}\frac{(u+g(s_{j-1,d}^{(0)}))^{k/2-1}}{q(u+g(s_{j-1,d}^{(0)})}e^{-( u+g(s_{j-1,d}^{(0)}) )/2}
%&\leq&M_1\varphi(u)\left(g(s_{j-1,d}^{(0)})\right)^{(k/2-1)_++1/\alpha+\eta}e^{-\frac{g(s_{j-1,d}^{(0)})}{2}}d.
\EQN as $u\to\IF.$
Next, for any $\alpha\in (0,2)$ or $\alpha=2, k>1$, since $\frac{t^{k/2-1}}{q(t)}$ is a regularly varying function at infinity with index $k/2-1+1/\alpha>0$, by Karamata's theorem (see, e.g., \cite{Res08}) we have % derive there exist constant $u_0>1$ such that for any $u>u_0$ and $j\geq N_{d,\delta}$,
\BQNY
\frac{(u+g(s_{j-1,d}^{(0)}))^{k/2-1}}{q(u+g(s_{j-1,d}^{(0)})}\leq  \mathcal{C }  \frac{(g(s_{j-1,d}^{(0)})^{k/2-1}}{q(g(s_{j-1,d}^{(0)})}\frac{u^{k/2-1}}{q(u)}
\EQNY
holds for all $u>u_0$ and all $j\geq N_{d,\delta}$,  with some  constant $u_0>1$.
 Moreover, we can write  $\frac{t^{k/2-1}}{q(t)}=t^{k/2-1+1/\alpha}\ell(t)$ with $\ell(t)$ %is a slowly varying function. Since $\ell\sim \ell^*$ as $t\rw\IF$ with $\ell^*$
 a normalized slowly varying function (see e.g., \cite{BI1989}). %, then we can assume that $\frac{t^{\frac{k}{2}-1}}{\overleftarrow{K}(t^{-1/2})}=t^{k/2-1+1/\alpha}\ell^*(t)$.
 One can check that $t^{k/2-1+1/\alpha}\ell(t)e^{-\frac{t}{2}}, t\ge t_0$ is decreasing for some large $t_0$.
Therefore, we  conclude that for $\delta$ and $d$ small enough  %and $u$ large enough
\BQN \label{eq:Idelt2}
\lim_{u\to\IF}\frac{I_{\delta,2}(u)}{\varphi(u)}&\leq&  \mathcal{C }  \sum_{j=N_{d,\delta}}^{\IF}\frac{(g(s_{j-1,d}^{(0)})^{k/2-1}}{q(g(s_{j-1,d}^{(0)})}e^{-\frac{g(s_{j-1,d}^{(0)})}{2}}
d\nonumber\\
&\leq&  \mathcal{C }  \sum_{j=N_{d,\delta}}^{\IF}\int_{-(j-1)d}^{-(j-2)d}\frac{\left(g(\invf(t))\right)^{k/2-1}}{q(g(\invf(t)))}e^{-\frac{g(\invf(t))}{2}}dt\\
%&=& M_2\int_{-\IF}^{f(\delta)+2d}\left(g(\invf(t))\right)^{(k/2-1)_++1/\alpha+\eta}e^{-\frac{g(\invf(t))}{2}}dt\\
&\le&  \mathcal{C }   \int_{0}^{\invf(f(\delta)+3d)}(C(t))^{1/\alpha}\frac{\left(g(t)\right)^{ k/2-1 }}{q(g(t))}e^{-\frac{g(t)}{2}}dt.\nonumber
\EQN
Consequently, letting $\delta\rw 0$, we derive by \textbf{C}(0) that  \eqref{eq:I23} holds for $I_{\delta,2}(u)$ when $\alpha\in (0,2)$ or $\alpha=2, k>1$.
For $\alpha=2, k=1$, we have from (\ref{K}) that
\BQNY
\frac{(u+g(s_{j-1,d}^{(0)})^{k/2-1}}{q(u+g(s_{j-1,d}^{(0)})}\leq  \mathcal{C }  _1,\ \ \frac{ u  ^{k/2-1}}{q(u)}\ge  \mathcal{C }  _2.
\EQNY
%with some positive constants $ \mathcal{C }  _1, \mathcal{C }  _2>0$.
Then, similar arguments as above show that \eqref{eq:I23} still holds for $I_{\delta,2}(u)$ when  $\alpha=2, k=1$.\\
Below we consider $I_{\delta,3}(u)$.
Let $Y(t)=X(\invf(t)),t\in[f(1-\delta), f(1))$. Then  correlation function of it is given by $r_Y(t,s)=r(\invf(t),\invf(s))$. It follows from \textbf{D}($1$) that for any $\delta$ small enough there exists some $M_3>1$ such that
\BQNY
\sup_{t\neq s\in (1-\delta,1)}\frac{1-r(t,s)}{K^2(|f(t)-f(s)|)}< M_3,
\EQNY
which means that %for $\delta$ small enough
\BQNY
1-r_Y(t,s)\leq M_3 K^2(|t-s|)\leq \frac{1}{2}K^2((8M_3)^{1/\alpha}|t-s|), \ \ \ s\neq t\in [f(1-\delta), f(1)).
\EQNY
%By the property of regularly varying function, we can rewrite the above inequality as
%\BQNY
%1-r_Y(s,t)\leq \frac{1}{2}K^2((8M_4)^{1/\alpha}|t-s|), \ \ \ s\neq t\in [f(T-\delta), f(T)).
%\EQNY
Furthermore, with the aid of (\ref{eqr}) we have  for $\delta$ small enough
\BQNY
1-r_Y( t,s)\leq 1-r((8M_3)^{1/\alpha}|t-s|),\ \ \ s\neq t\in [f(1-\delta), f(1)).
\EQNY
Thus, similarly as the above case for $I_{\delta,2}(u)$, it follows that %from Lemma \ref{LEV} and Lemma \ref{LES} that for $\delta$ small enough and $u$ large enough,
\BQN   \label{eqt}
I_{\delta,3}(u)&\leq& % 2^n\mathbb{P}\left(\sup_{t\in(f(0),f(\delta)]}Y_{\vk{b}}^2(t)>u\right)
 2^n\mathbb{P}\left(\sup_{t\in[f(1-\delta), f(1)) }V_{\vk{b}}^2((8M_3)^{1/\alpha}t)>u\right)\nonumber\\
&\leq&2^n\mathbb{P}\left(\sup_{t\in[0,( f(1)-f(1-\delta))(8M_3)^{1/\alpha}]}V_{\vk{b}}^2(t)>u\right)\\
&\leq& 2^{n+1}( f(1)-f(1-\delta))(8M_3)^{1/\alpha}\mathcal{H}_{\alpha}\frac{2^{1-k/2}}{\Gamma(k/2)}\varphi(u)\nonumber
\EQN
holds for $\delta$ small enough and $u$ large enough. % where $Y_{\vk{b}}(t)=\sum_{i=1}^nb_i^2Y_i^2(t), t\in\mathbb{R}$ with $Y_i, 1\leq i\leq n$ independent copies of $Y$ and $1=b_1=\cdots=b_k> b_{k+1} \cdots\geq b_n>0$.
Consequently, in light of $f(1)<\IF$ we derive  that  \eqref{eq:I23} holds for $I_{\delta,3}(u)$.
The   proof  is complete.\QED

\proofkorr{cor} By direct calculation, we have
\BQNY
1-\E{\chi(t)\chi(s)})&=&\sqrt{\frac{s}{t}\frac{1-t}{1-s}}\left(\sqrt{(1+\frac{t-s}{s})(1+\frac{t-s}{1-t})}-1\right)\\
&=&\frac{|t-s|}{2t_0(1-t_0)}(1+o(1)),\ \ \ s<t\rw t_0\in(0,1),
\EQNY
which, in the notation of \eqref{INCO}, means that $C(t)=\frac{1}{2t(1-t)}$, $K^2(t)=t$ and $\alpha=1$.
Further
\BQNY
%\int_0^1\frac{1}{2t(1-t)}dt=\IF,\ \ \
f(t)=\int_{1/2}^t\frac{1}{2s(1-s)}ds=\frac{1}{2}\ln\frac{t}{1-t},\ \ f(1)=-f(0)=\IF.
\EQNY
Thus, we have
\BQNY
s_{j,d}^{(0)}=\invf(-jd)=\frac{1}{1+e^{2jd}},\ \ s_{j,d}^{(1)}=\invf(jd)=\frac{1}{1+e^{-2jd}} \ \ j\in \mathbb{N}\cup\{0\}, d>0.
\EQNY
Next we verify  conditions of Theorem \ref{TH1} and Theorem \ref{Addition}. Let therefore in the following  $S\in\{0,1\}$.
%Clearly, condition {\bf A}$(S)$ is satisfied since $g_\nu(t) \uparrow \IF$ as $t\to0$ and $t\to 1$.
Elementary calculations yield that
\BQNY
g_\nu'(t)=\frac{1+c^2(t)+2\nu c(t)}{1+c^2(t)}\frac{2t-1}{t(1-t)(1-\ln(4t(1-t)))}.
\EQNY
Since $1-\ln(4t(1-t))\geq 1$ and $c(t)\geq 0$ for $t\in(0,1)$, we have $g_\nu'(t)<0, t\in(0,1/2)$ and $g_\nu'(t)>0, t\in(1/2,1)$, implying that {\bf A}$(S)$ is satisfied since $g_\nu(t) \uparrow \IF$ as $t\to0$ and $t\to 1$.
Now, for any $s<t$ and $s,t\in \Delta_{j,d}^{(0)}=[\frac{1}{1+e^{2jd}},\frac{1}{1+e^{2(j-1)d}}]$ we have
\BQNY
 e^{-2d}\leq\frac{s}{t}\leq 1,\ \ \ \frac{1}{2}\leq\frac{1-t}{1-s}\leq 1,\ \ \ \frac{t-s}{s}\leq e^{2d}-1,\ \ \ \frac{t-s}{1-t}\leq e^{2d}-1.
\EQNY
Thus, there exists $d_0>0$ such that for $s<t\in \Delta_{j,d_0}^{(0)}=[\frac{1}{1+e^{2jd_0}},\frac{1}{1+e^{2(j-1)d_0}}]$
\BQNY
\frac{1-r(s,t)}{|f(t)-f(s)|}&=&\frac{2\sqrt{\frac{s}{t}\frac{1-t}{1-s}}\left(\sqrt{(1+\frac{t-s}{s})(1+\frac{t-s}{1-t})}-1\right)}
{\ln(1+\frac{t-s}{s})+\ln(1+\frac{t-s}{1-t})}\\
&\leq& 4\frac{\frac{t-s}{s}+\frac{t-s}{1-t}}{\frac{1}{2}(\frac{t-s}{s}+\frac{t-s}{1-t})}=8
\EQNY
and
\BQNY
\frac{1-r(s,t)}{|f(t)-f(s)|}\geq e^{-d_0}\frac{\frac{t-s}{s}+\frac{t-s}{1-t}}{2(\frac{t-s}{s}+\frac{t-s}{1-t})}\geq \frac{e^{-d_0}}{2}.
\EQNY
Moreover, for all $t\in \Delta_{j,d_0}^{(0)}, s\in \Delta_{j+l,d_0}^{(0)}$, with $l\in \mathbb{N}$ large enough,
$$
|r(s,t)|=\frac{s(1-t)}{\sqrt{st(1-t)(1-s)}}\leq 2\sqrt{\frac{s}{t}}\leq 2e^{-(l-1)d_0}\leq  \mathcal{C }   l^{-2}.
$$
Similarly, we can also derive that for $s<t\in \Delta_{j,d_0}^{(1)}=[\frac{1}{1+e^{-2(j-1)d_0}},\frac{1}{1+e^{-2jd_0}}]$
\BQNY
\frac{e^{-d_0}}{2}\leq \frac{1-r(s,t)}{|f(t)-f(s)|}\leq 8,
\EQNY
and
for all $t\in \Delta_{j,d_0}^{(1)}, s\in \Delta_{j+l,d_0}^{(1)}$, with $l \in \mathbb{N}$ large enough,
$$
\abs{r(s,t)}\leq  \mathcal{C }  l ^{-2}.
$$
Therefore, conditions {\bf B}(S), {\bf E}(S) are satisfied.
%Denote $a(t)\sim b(t)$ if $\lim_{t\rw 0}\frac{a(t)}{b(t)}=1$ and by d
Direct calculations show that
\BQNY
c(t)\sim \ln\ln\frac{1}{t},\ \ \ g_\nu(t)\sim\ln\ln\frac{1}{t},\ \ \ e^{-g_\nu(t)}\sim\frac{1}{(\ln\ln\frac{1}{t})^{2\nu}\ln\frac{1}{t}}
\EQNY
as $t\to 0.$
Consequently,
\BQNY
\frac{(g_\nu(t))^{1/2}}{t(1-t)}e^{-g_\nu(t)}\sim\frac{1}{t\ln\frac{1}{t}(\ln\ln\frac{1}{t})^{2\nu-1/2}}
\EQNY
as $t\to 0.$
This implies that $\int_{1/2}^{1}\frac{(g_\nu(t))^{1/2}}{t(1-t)}e^{-g_\nu(t)}dt=\int_0^{1/2}\frac{(g_\nu(t))^{1/2}}{t(1-t)}e^{-g_\nu(t)}dt<\IF$ if and only if $\nu>3/4$. Consequently, the claim follows by an application of \netheo{TH1} and Theorem \ref{Addition}.
 This completes the proof.
 \QED

%\proofkorr{cor4} Thanks to \nelem{A2} and \netheo{TH1}, we   only need to show that   condition  {\bf B}(0) is satisfied for the correlation function of Brownian motion.  This will be verified as  a special case (with $H=1/2$)  in the proof of Corollary \ref{cor2} below.\QED

\proofkorr{cor2}  As discussed in Introduction,
the correlation of the fBm satisfies
\BQNY
1-\text{Corr}(B_H(t),B_H(s))= \frac{|t-s|^{2H}}{2t_0^{2H}}(1+o(1)), \ \ s,t\rw t_0\in(0,1],
\EQNY
which means that
\BQNY
C(t)=\frac{1}{2 t^{2H}},\ \ K^2(t)= t^{2H},\ \ \alpha =2H.
\EQNY
%Direct calculation yields that
%\BQNY
%r(s,t)=Corr(B(t),B(s))=1-\frac{|t-s|}{\sqrt{t\vee s}(\sqrt{t}+\sqrt{s})}=1-\frac{|t-s|}{2t_0}(1+o(1)),\ \ s,t\rw t_0 \in (0,1],
%\EQNY
%where $t\vee s=\max(s,t)$.
We only need to verify the conditions {\bf B}(0) and {\bf E}(0). It follows that
\BQNY
&& \ \ f(t)=2^{-\frac{1}{2H}}\ln(2t),\ \ s_{j,d}^{(0)}=\invf(-jd)=\frac{1}{2}e^{-2^{\frac{1}{2H}}jd},\\ &&\Delta_{j,d}^{(0)}=[\frac{1}{2}e^{-2^{\frac{1}{2H}}jd},  \frac{1}{2}e^{-2^{\frac{1}{2H}}(j-1)d}], \ j\in \mathbb{N}.
\EQNY
%We concentrate on the interval $(0,1/2]$ since for the interval $[1/2,1]$ it's the standard chi-square with trend problem.
Without loss of generality, hereafter we assume that $s<t$. For any $s,t\in \Delta_{j,d}^{(0)}, j\in \mathbb{N}$ we have $\frac{ t-s }{s}\leq e^{2^{\frac{1}{2H}}d}-1$, which implies that there exists $d_0>0$ such that %for  any $0<d<d_0$ %, we have
\BQNY
 \frac{ t-s }{2s}\leq  |\ln (1+\frac{t-s}{s})|\leq \frac{2(t-s) }{s}
 \EQNY
 holds for all $s,t\in \Delta_{j,d_0}^{(0)}, j\in \mathbb{N}$. Moreover, for $s,t\in \Delta_{j,d_0}^{(0)}, j\in \mathbb{N}$ \HH{with $d_0$ small enough}, %$0<d<d_0$,
 $$
 \frac{(t^H-s^H)^2}{|t-s|^{2H}}=\frac{\left(1-(\frac{s}{t})^{H}\right)^2}{|\frac{t-s}{t}|^{2H}}\leq 2H^2\left(\frac{t-s}{t}\right)^{2-2H}<1/2.$$
Thus, we have, for %any $0<d<d_0$ and
all $s,t\in \Delta_{j,d_0}^{(0)}, j\in \mathbb{N}$ \HH{with $d_0$ small enough}
\BQNY
\frac{1-r(s,t)}{|f(t)-f(s)|^{2H}}\leq\frac{|t-s|^{2H}}{t^{H}s^H|\ln (1+\frac{t-s}{s})|^{2H}}\leq 2^{2H}\frac{s^H}{t^H}\leq 2^{2H},
\EQNY
and
$$\frac{1-r(s,t)}{|f(t)-f(s)|^{2H}}\geq\frac{|t-s|^{2H}}{2t^{H}s^H|\ln (1+\frac{t-s}{s})|^{2H}}\geq 2^{-2H-1}\frac{s^H}{t^H} \geq 2^{-2H-1}e^{-2^{\frac{1}{2H}}d_0H}>0.$$
In addition, %for  $0<d<d_0$ and $l\geq l_0\in \mathbb{N}$ with $l_0$ sufficiently large,
\BQNY
|r(s,t)|&=&\frac{|t^{2H}+s^{2H}-|t-s|^{2H}|}{t^Hs^H}=\frac{|(\frac{s}{t})^{2H}+1-(1-\frac{s}{t})^{2H}|}{(\frac{s}{t})^H}\leq \frac{(\frac{s}{t})^{2H}+4H\frac{s}{t}}{(\frac{s}{t})^H} \\
&\leq&  \mathcal{C }  \left(e^{-2^{\frac{1}{2H}}dH l}+e^{-2^{\frac{1}{2H}}d(1-H)l}\right)\leq  \mathcal{C }  l^{-2},
\EQNY
holds for  $t\in \Delta_{j,d_0}^{(0)}, s\in \Delta_{j+l,d_0}^{(0)}, j\in \mathbb{N}$, with $l$ sufficiently large and \HH{with $d_0$ small enough.}
 Consequently,  both  conditions {\bf B}(0) and {\bf E}(0) are satisfied, and thus % Further, condition {\bf C}(0) holds since $\int_{0}^1 \frac{1}{t}(g(t))^{k/2+1/(2H)- 1}e^{-\frac{g(t)}{2}}dt<\IF$.
the claim follows  by applying Theorem \ref{TH1} and Theorem \ref{Addition}.
 %with $\eta=0$, we can derive (\ref{e1}). Thus we
 This completes the proof.
 \QED

\prooftheo{Thd0} The proof   is  similar to that of Theorem \ref{Th0}. The main difference is the expansion of the correlation function of $Y_{\vk{1}}(t,\theta)$ defined in (\ref{Y}). Here it is given by (compare with (\ref{PRC2}))
\BQNY
\Corr\left(Y_{\vk{1}}(t,\vk{\theta}),Y_{\vk{1}}(s,\vk{\theta}^{'})\right)&=&1-\left(C_1(t_0)\prod_{j=2}^n\cos^2\theta_j+\sum_{i=2}^kC_i(t_0)\sin^2\theta_i
\prod_{j={i+1}}^n\cos^2\theta_j\right)K^2(|t-s|)(1+o(1))\\
&&-\sum_{i=2}^{n-1}\frac{1}{2}\left(\prod_{l=i+1}^n\cos^2 \theta_l\right)(\theta_i-\theta_i^{'})^2(1+o(1))
-\frac{1}{2}(\theta_n-\theta_n^{'})^2(1+o(1))
\EQNY
as $t,s\rw t_0, |\theta_i-\theta_i^{'}|\rw 0, 2\leq i\leq n.$
The rest of the proof is the same as that in the proof of Theorem \ref{Th0}. This completes the proof. \QED

\prooftheo{TH2} The proof follows with the same  arguments as those in the proof  of Theorem \ref{TH1}, and thus being omitted here.\QED

\proofkorr{cor3} %Assume that for $s,t\rw t_0\in (0,1)$
Denote for any $s,t\in(0,1)$
\BQNY
 r_1(s,t)=\text{Corr}(B(t),B(s)),
 r_2(s,t)=\text{Corr}(W(t),W(s)),  r_3(s,t)=\text{Corr}(B_H(t),B_H(s)).
\EQNY
In view of the proof of Corollary \ref{cor} and Corollary \ref{cor2}, we have (in the notation of Section 3)
\BQNY
 C_1(t)=\frac{1}{2t(1-t)},\ C_2(t)=\frac{1}{2t},\ C_3(t)=\frac{1}{2^{2H}t^{2H}},\
  K^2(t)=t,
 \ K_3^2(t)=2^{2H-1}t^{2H},\ k=2<3=n.
\EQNY
Then, we have for $H\in(1/2,1)$
\BQNY
&&C^*(t)=\max\left(\frac{1}{2t(1-t)}, \frac{1}{2t}, \frac{1}{2t}\right)=\frac{1}{2t(1-t)}, \   f^*(t)=\frac{1}{2}\ln\frac{t}{1-t},\\
&& s_{j,d}^{(0)}=\invfs(-jd)=\frac{1}{1+e^{2jd}},\ s_{j,d}^{(1)}=\invfs(jd)=\frac{1}{1+e^{-2jd}} \ \ j\in \mathbb{N}\cup\{0\}, d>0.
\EQNY
Without loss of generality, we assume that $s<t$.
Since $\frac{t-s}{s}\leq e^{2d}-1$ and $\frac{t-s}{1-t}\leq e^{2d}-1$ hold for $s,t\in \Delta_{j,d}^{(0)}$ or $s,t\in \Delta_{j,d}^{(1)}, \ \ j\in N$, we have there exists some $d_0>0$ such that % for any $0<\delta<\delta_0$
\BQNY
\ln(1+\frac{t-s}{s})\geq \frac{t-s}{2s},\ \ \ln(1+\frac{t-s}{1-t})\geq \frac{t-s}{2(1-t)},\   1/2<\frac{s}{t} \leq 1
\EQNY
hold for $s,t\in \Delta_{j,d_0}^{(0)}$ or $s,t\in \Delta_{j,d_0}^{(1)}, \ \ j\in \mathbb{N}$. This implies,  for any $s<t\in \Delta_{j,d_0}^{(0)}$ or any $s<t\in \Delta_{j,d_0}^{(1)}$    %which lead to, with the same range of $s,t$,
\BQNY
 \frac{1-r_2(s,t)}{|f^*(t)-f^*(s)|}=\frac{2(t-s)}{\sqrt{t}(\sqrt{t}+\sqrt{s})|\ln(1+\frac{t-s}{s})+\ln(1+\frac{t-s}{1-t})|}\leq \frac{4s(1-t)}{\sqrt{t}(\sqrt{t}+\sqrt{s})(1-t+s)}\leq 2
 \EQNY
 and
 \BQNY
&&\frac{1-r_3(s,t)}{|f^*(t)-f^*(s)|}=\frac{|t-s|^{2H}-(t^H-s^H)^2}{t^Hs^H|\ln(1+\frac{t-s}{s})+\ln(1+\frac{t-s}{1-t})|}\leq \frac{|t-s|^{2H}}{t^Hs^H|\ln(1+\frac{t-s}{s})+\ln(1+\frac{t-s}{1-t})|}\\
&&\leq \frac{2|t-s|^{2H-1}s(1-t)}{t^Hs^H(1-t+s)}\leq  \frac{2|t-s|^{2H-1}s^{1-H}}{t^H}\leq\sup_{x\in(1/2,1]}2|1-x|^{2H-1}x^{H-1}<\IF.
\EQNY
Thus, condition {\bf B}$(S)$ with $S=0, 1$ holds for $r_2(\cdot,\cdot),r_3(\cdot,\cdot)$. Additionally, it has been shown in proof of Corollary \ref{cor} that condition {\bf B}$(S)$ with $S=0, 1$ holds for $r_1(\cdot,\cdot)$. Consequently, in the light of Theorem \ref{TH2}  we derive that, as $u\to\IF$
\BQNY
&&\pc{\sup_{t\in(0,1)}\left(\frac{B^2(t)}{t(1-t)}+\frac{W^2(t)}{t}+\frac{B_H^2(t)}{t^{2H}}-g(t)\right)>u}\\
&&\sim (2\pi)^{-3/2} u^{3/2}e^{-u/2}\int_0^1dt\int_{-\pi}^\pi d\theta_2\\
&&\times\int_{-\pi/2}^{\pi/2}\left(\frac{1}{2t(1-t)}\cos^2(\theta_2)\cos^2(\theta_3)+\frac{1}{2t}
\sin^2(\theta_2)\cos^2(\theta_3)\right)\cos(\theta_3)e^{-\frac{g(t)}{2}}d\theta_3\\
&&= \frac{u^{3/2}e^{-u/2}}{3\sqrt{2\pi}}\int_0^1\frac{2-t}{t(1-t)}e^{-\frac{g(t)}{2}}dt.
\EQNY
This completes the proof. \QED

\section{Appendix}

This section is concerned with the proof of \netheo{Addition}. We first present, in the following theorem, an
 asymptotic 0-1  behavior of the chi-square process $\chi_{\vk{b}}^2$, which is complementary to the results for Gaussian processes discussed in \cite{watanabe1971,QueWat72, Wat70}. \aH{Moreover, it can also be viewed as a generalization of the Kolmogorov-Dvoretsky-Erd\H{o}s integral test theorem; see, e.g., page 163 in \cite{IM1974}  or Theorem A in \cite{KLS96}. Recall that, for the  Bessel process  $||\vk{W}_n(t)||:=\sqrt{\sum_{i=1}^n{W_i^2(t)}}, t\geq 0$   defined in Corollary \ref{cor4}, the
Kolmogorov-Dvoretsky-Erd\H{o}s integral test theorem tells that for any positive continuous ultimately increasing $g(\cdot)$ when $t\rw 0$,
$$\pk{ ||\vk{W}_n(t)||\leq  \sqrt{tg(t)}, \ \ \text{ultimately as}\ \ t\rw 0}=1\quad \text{or}\quad 0$$
holds
according as $$\int_0^1\frac{1}{t}(g(t))^{n/2}e^{-\frac{g(t)}{2}}dt<\IF \quad \text{or} \ =\IF.$$
 %which is a specific corollary of  the following theorem with $\chi_{\vk{b}}^2(t)=\vk{W}_n^2(t)$ followed by the proof of Corollary \ref{cor2}.
 }
\BT\label{THM} Let $\{X(t),t\in E\}$ be given as in  \netheo{Th0} with $E=(0,1)$, and $S\in\{0,1\}$. Then,
%Suppose that \textbf{Condition A}($S$) and \textbf{Condition E}($S$) are satisfied.\\
for any nonnegative continuous function $g(t)$ satisfying {\bf A}(S),
\BQN\label{eq:chig1}
\pk{\chi_{\vk{b}}^2(t)\leq g(t) \ \ \text{ultimately as}\ \ t\rw S }=1
\EQN
holds provided that $I_g(S)<\IF$, $\abs{f(S)}=\IF$ and condition   {\bf B}(S) is satisfied, or $\abs{f(S)}<\IF$ %(implying $I_g(S)<\IF$)
and condition {\bf D}(S) is satisfied, and
\BQN\label{eq:chig2}
\pk{\chi_{\vk{b}}^2(t)\leq g(t) \ \ \text{ultimately as}\ \ t\rw S }=0
\EQN
holds provided that
 $I_g(S)=\IF$ and condition {\bf E}(S) is satisfied.
\ET

\textbf{Proof.} The idea of the proof comes from \cite{watanabe1971, QueWat72, Wat70}. Without loss of generality, we present only the proof for the case $S=0$. We prove first that \eqref{eq:chig1} holds if $I_g(0)<\IF$, $f(0)=-\IF$ and condition   {\bf B}(0) is satisfied.
Let
$$E_{j,d}=\left\{\sup_{t\in \Delta^{(0)}_{j,d}}\chi_{\vk{b}}^2(t)> g(s_{j-1,d}^{(0)}) \right\}, \ \  j\in \mathbb{N}, d>0.$$
%We begin with the scenario that $I_g(S)<\IF$.
%Note that $\frac{t^{\frac{k}{2}-1}}{\overleftarrow{K}(t^{-1/2})}$ is a regularly varying function at $\IF$ with index $k/2-1+1/\alpha$ and thus we can write it as $t^{k/2-1+1/\alpha}\ell(t)$ with $\ell(t)$ is a slowly varying function. Since $\ell\sim \ell^*$ as $t\rw\IF$ with $\ell^*$ a normalized slowly varying function (see e.g., \cite{BI1989}), then we can assume that $\frac{t^{\frac{k}{2}-1}}{\overleftarrow{K}(t^{-1/2})}=t^{k/2-1+1/\alpha}\ell^*(t)$. One can check that $t^{k/2-1+1/\alpha}\ell^*(t)e^{-\frac{t}{2}}$ is decreasing for $t$ large enough.
%Thus it follows from (\ref{eqr})-(\ref{eqs}) that
Similarly to the derivation of \eqref{eq:Idelt2}, we have, for $\delta, d$ small,
\BQN\label{con1}
\sum_{j=N_{d,\delta}}^\IF\mathbb{P}\left( E_{j,d}\right)&\leq& 2^n\sum_{j=N_{d,\delta}}^\IF
 H_\alpha\frac{2^{1-k/2}}{\Gamma(k/2)}(4M_1)^{1/\alpha}d\ \varphi(g(s_{j-1,d}^{(0)}))\nonumber\\
 &\leq& \mathcal{C }   \sum_{j=N_{d,\delta}}^\IF \frac{(g(s_{j-1,d}^{(0)}))^{\frac{k}{2}-1}}{ q(g(s_{j-1,d}^{(0)}))  }e^{-\frac{g(s_{j-1,d}^{(0)})}{2}}d\nonumber\\
 &\leq& \mathcal{C }   \sum_{j=N_{d,\delta}}^{\IF}\int_{-(j-1)d}^{-(j-2)d}\frac{(g(\invf(t)))^{k/2-1}}{q (g(\invf(t))) }e^{-\frac{g(\invf(t))}{2}}dt\nonumber\\
 &\leq&  \mathcal{C }   \int_{0}^{1/2} (C(t))^{1/\alpha}\frac{(g(t))^{\frac{k}{2}-1}}{q(g(t)) }e^{-\frac{g(t)}{2}}dt<\IF
\EQN
holds when $\alpha\in(0,2)$ or $\alpha=2, k>1$, and
 \BQNY
\sum_{j=N_{d,\delta}}^\IF\mathbb{P}\left( E_{j,d}\right) \leq     \mathcal{C }  _1 \int_{0}^{1/2} (C(t))^{1/\alpha}e^{-\frac{g(t)}{2}}dt<\IF
\EQNY
holds when $\alpha=2, k=1$. Thus, by Borel-Cantelli lemma
$$
\pk{\exists j_g:  \sup_{t\in \Delta^{(0)}_{j,d}}\chi_{\vk{b}}^2(t)\le g(s_{j-1,d}^{(0)})\ \mathrm{for\ all\ } j\ge j_g }=1
$$
which implies \eqref{eq:chig1} since $g(t)\uparrow$  as $t\downarrow 0$. Next we  prove \eqref{eq:chig1} under the conditions that  $f(0)>-\IF$ and {\bf D}(0) is satisfied. Let for any fixed small $\delta>0$, $t_n= \delta/n, n\in \NN $. Denote
$$
H_{n,\delta}=\left\{\sup_{t\in[t_{n+1},t_n]} \chi_{\vk{b}}^2(t)> g(t_n)  \right\}.
$$
Similar arguments as in \eqref{eqt} yield that, for $N_0$ sufficiently large
\BQNY
\sum_{n=N_0}^\IF \pk{H_{n,\delta}}&=& \sum_{n=N_0}^\IF \pk{ \sup_{t\in[f(t_{n+1}),f(t_n)]} \chi_{\vk{b}}^2(\overleftarrow{f}(t))> g(t_n)  }\\
&\le &   \mathcal{C }    \sum_{n=N_0}^\IF (f(t_n)-f(t_{n+1})  ) \varphi(g(t_n))\\
&\le&   \mathcal{C }  _1  (f(t_{N_0})-f(0)  )<\IF
\EQNY
holds. %, with some positive constants $ \mathcal{C }  , \mathcal{C }  _1$.
Thus, the proof of \eqref{eq:chig1}  is again completed by Borel-Cantelli lemma.
Finally, we show  \eqref{eq:chig2}  when $I_g(0)=\IF$ and condition {\bf E}(0) is satisfied. Note that
\BQN \label{eq:chibg}
&&\pk{\chi_{\vk{b}}^2(t)\leq g(t) \ \ \text{ultimately as}\ \ t\rw 0 }\nonumber\\
&&\le \pk{\sup_{\vk{v}\in\mathcal{S}_{k-1} }Y_k(t,\vk{v})\leq \sqrt{g(t)} \ \ \text{ultimately as}\ \ t\rw 0 }\\
&&=\pk{\sup_{\vk{v}\in\mathcal{S}_{k-1} }Y_k(\overleftarrow{f}(-t),\vk{v})\leq \sqrt{g(\overleftarrow{f}(-t))} \ \ \text{ultimately as}\ \ t\rw \IF },\nonumber
\EQN
%The claim is equivalent to
%$$\pk{\sup_{\vk{v}\in\mathcal{S}_{k-1} }Y_k(t,\vk{v})\leq \sqrt{g(t)} \ \ \text{ultimately as}\ \ t\rw 0 }=0,$$
where $ Y_k(t,\vk{v})=\sum_{i=1}^k X_i(t)v_i, \ \ (t,\vk{v})\in (0,1/2]\times\mathcal{S}_{k-1},$
with %$X_i, 1\leq i\leq k$ being independent copies of $X$ and
$\mathcal{S}_{k-1}$  the $(k-1)$-dimensional unit sphere.
By similar arguments as in Lemma 1.4 in \cite{watanabe1971} or Lemma 3.1  in \cite{QueWat72}, one can for free assume that
\BQN\label{Restriction}
2\ln t\leq g(\overleftarrow{f}(-t))\leq 3\ln t,\ \ \  t\in (1, \IF).
\EQN
Denote
\BQNY
F_{j,d}&=& \left\{ (j-1)d+ l\ q(g(s_{j,d}^{(0)})),\ \  l=0,1,\cdots, \left\lfloor\frac{d}{q(g(s_{j,d}^{(0)}))}\right\rfloor \right\},\\
G_{j,d}&=&   \left\{\left(l_1h_{j,d}, l_2h_{j,d}, \dots, l_{k-1}h_{j,d}, \pm\sqrt{1-\sum_{i=1}^{k-1}l_i^2h_{j,d}^2} \right ), \ \ l_1,\cdots, l_{k-1}\in \mathbb{Z}\right\},
\EQNY
with $h_{j,d}=1/\sqrt{g(s_{j,d}^{(0)})}$. Further, define
\BQNY
W_{j,d}=\left\{\sup_{ t \in F_{j,d} } \sup_{ \vk{v} \in G_{j,d}\cap  \mathcal{S}_{k-1}}Y_k( \overleftarrow{f}(-t),\vk{v})\le \sqrt{g(s_{j,d}^{(0)})}\right\}.
\EQNY
In view of \eqref{eq:chibg}, to prove \eqref{eq:chig2} it is sufficient that
%$$G_{j,d}=\left(\Delta_{j,d}^{(0)}\times \mathcal{S}_{k-1}\right)\bigcap F_{j,d}, \ \  \ \
%W_{j,d}=\left\{\sup_{(t,\vk{v})\in G_{j,d}}Y(t,\vk{v})<\sqrt{g(s_{j,d}^{(0)})}\right\}$$ $$ W_{j,d}^c=\left\{\sup_{(t,\vk{v})\in G_{j,d}}Y(t,\vk{v})\geq \sqrt{g(s_{j,d}^{(0)})}\right\}, \ \ \text{and} \ \ \mathcal{S}_{k-1}^{\epsilon}= \mathcal{S}_{k-1}\bigcap ([-\epsilon, \epsilon]^{k-1}\times [0,1]).$$
%It suffice to prove that there exists $k_0\in \mathbb{N}$ such that
$$
\pk{  W_{il_0,d}^c \ \ \mathrm{i.o.\ }}=1
$$
holds, with $l_0$ as in condition {\bf E}. It is noted that
\BQN\label{A3}
1-\pk{ W_{il_0,d}^c \ \ \mathrm{i.o.\ } }=\lim_{l\rw\IF}\prod_{i=l}^\IF\pk{ W_{il_0,d}}+\lim_{l\rw\IF}\left(\pk{\bigcap_{i=l}^\IF W_{il_0,d}}-\prod_{i=l}^\IF\pk{ W_{il_0,d}}\right).
\EQN
It remains to show that both limits above are equal to 0. For the first limit to be zero it is sufficient to show that $\sum_{i=1}^\IF \pk{ W_{il_0,d}^c }=\IF$. For any small $\eta>0$ denote $\mathcal{S}_{k-1}^{\eta}= \mathcal{S}_{k-1}\cap ([-\eta, \eta]^{k-1}\times [0,1])$.
Clearly,
$$
\pk{  W_{il_0,d}^c }\ge \pk{ \sup_{ t \in F_{j,d} } \sup_{ \vk{v} \in G_{j,d}\cap  \mathcal{S}^\eta_{k-1}}Y_k( \overleftarrow{f}(-t),\vk{v})> \sqrt{g(s_{j,d}^{(0)})}  }.
$$
We have from condition \textbf{E}(0)  that there exist some  $M_4>0, \delta_3>0$ and $d_2>0$ such that for any $\delta\in(0,\delta_3)$, $d\in(0,d_2)$,
\BQNY
1-r(\finv(-t),\finv(-s))\geq M_4 K^2(|t-s|)\geq 2K^2((4^{-1}M_4)^{1/\alpha}|t-s|), ~~t,s\in [(j-1)d, jd]
\EQNY
holds for all $j\geq N_{d,\delta}$. Let $r_{Y_k}(\cdot,\cdot,\cdot, \cdot)$ be the covariance function of $Y_k(\overleftarrow{f}(-t),\vk v)$. By using similar arguments as in the proof of Lemma 3.3 in \cite{HAJI2014} we have, for sufficiently small $\eta$
%%%%%%%%%%%%%%%%%555
\COM{

%Let $r(\cdot)$ be a covariance function of a stationary Gaussian process such that for some constant $ \delta_2>0$
% \BQNY
%  1-r(t)\leq 2K^2(\abs{t}),\ \ \ t\in(0,\delta_2).
% \EQNY
% See e.g.,  \nelem{LEV}.
Therefore, for any   $d\in(0,d_1)$ and $\delta\in(0,\min(\delta_1,\delta_2))$
\BQNY
%\E{\left(X(t)\right)^2}&=&\E{\left(V((4M_1)^{1/\alpha}f(t))\right)^2},\\
\mathbb{E}\left(X(t)X(s)\right) \leq  \mathbb{E}\left(V((4^{-1}M_1)^{1/\alpha}f(t))V(((4^{-1}M_1)^{1/\alpha}f(s))\right)
\EQNY
holds for all $t, s\in\Delta_{j,d}^{(0)}$ and $j\geq N_{d,\delta}$, where $V$ is the stationary Gaussian process given in \nelem{LEV}.
%Apparently, we can derive the desired correlation inequality by the above inequality and (\ref{eqr}).
Therefore, \nelem{LES} in Appendix  leads to
\BQNY
\mathbb{P}\left(W_{j,d}^c\right)&\geq& 2^{-k} \mathbb{P}\left(\sup_{(t,\vk{v})\in G_{j,d}}Z(4^{-1}M_1)^{1/\alpha}f(t),\vk{v})>g(s_{j,d}^{(0)})\right)\\
&=&2^{-k} \mathbb{P}\left(\sup_{(t,\vk{v})\in G_{j,d}^*}Z(t,\vk{v})>g(s_{j,d}^{(0)})\right)
\EQNY
holds for $ \delta\in(0,\delta_3)$, $ d\in(0,d_2)$ and $j\geq N_{d,\delta}$, where
$$Z(t,\vk{v})=\sum_{i=1}^k V_i(t)v_i, \ \ (t,\vk{v})\in (0,1/2]\times\mathcal{S}_{k-1} $$
with $V_{i}$ being independent copies of the Gaussian process in \nelem{LEV} and
$$ G_{j,d}^*=\left([0, (4^{-1}M_1)^{1/\alpha}d]\times \mathcal{S}_{k-1}^\epsilon\right)\bigcap F_{j,d}^*$$
with
\BQNY
F_{j,d}^*&=& \left\{l(4^{-1}M_1)^{1/\alpha}(\ln g(s_{j,d}^{(0)}))^{-1} \overleftarrow{K}((g(s_{j,d}^{(0)}))^{-1/2}), l\in \mathbb{Z}\right\}\\
&& \ \ \times \left\{\left(l_1h_{j,d}, l_2h_{j,d}, \dots, l_{k-1}h_{j,d}, \pm\sqrt{1-\sum_{i=1}^{k-1}l_i^2h_{j,d}^2} \right ), \ \ l_1,\dots, l_{k-1}\in \mathbb{Z}\right\}.
\EQNY
%Further$\Delta_d=[0,(4M_1)^{1/\alpha}d]$ and
%\BQNY
%I_{\delta,2}(u)\leq 2^n\sum_{j=N_{d,\delta}}^{\IF}\mathbb{P}\left(\sup_{t\in \Delta_{d}}V^2_{\vk{b}}(t)>u+g(s_{j,d}^{(0)})\right).
%\EQNY
%With the aid of
}
%%%%%%%%%%%%%%%%%%%%%%%%%%%%%%%%%%%%%%%%%%%%%%%%%%%%%%%%%%%%%%%%%%%%%%%%%%%%%%%%%%%%%55
\BQNY
1- r_{Y_k}(s,\vk{v}, t,\vk{v}') \ge 2K^2((4^{-1}M_4)^{1/\alpha}|t-s|) +\frac{3}{8} \sum_{i=1}^{k-1}(v_i-v_i')^2, \ \ t,s\in [(j-1)d, jd], \vk{v}, \vk{v}'\in  \mathcal{S}^\eta_{k-1}
\EQNY
holds for all  $j\geq N_{d,\delta}$, and in the light of Slepain's inequality
\BQNY
\pk{ \sup_{ t \in F_{j,d} } \sup_{ \vk{v} \in G_{j,d}\cap  \mathcal{S}^\eta_{k-1}}Y_k( \overleftarrow{f}(-t),\vk{v})> \sqrt{g(s_{j,d}^{(0)})}  }\ge
\pk{ \sup_{ t \in \widetilde{F}_{j,d} } \sup_{ \widetilde{\vk{v}} \in \widetilde{G}_{j,d}\cap [-\eta,\eta]^{k-1}}Z((4^{-1}M_1)^{1/\alpha} t,\widetilde{\vk{v}})> \sqrt{g(s_{j,d}^{(0)})}  },
\EQNY
where $\widetilde{\vk{v}}=(v_1,\dots, v_{k-1})\in\R^{k-1}$,
\BQNY
\widetilde{F}_{j,d}&=& \left\{  l\ q(g(s_{j,d}^{(0)})),\ \  l=0,1,\cdots, \left\lfloor\frac{d}{q(g(s_{j,d}^{(0)}))}\right\rfloor \right\},\\
\widetilde{G}_{j,d}&=&   \left\{\left(l_1h_{j,d}, l_2h_{j,d}, \dots, l_{k-1}h_{j,d} \right ), \ \ l_1,\cdots, l_{k-1}\in \mathbb{Z}\right\},
\EQNY
and $Z$ is a centered homogeneous Gaussian random field with a.s. continuous sample paths, variance 1 and correlation function such that
\BQN\label{eq:corZ}
1-Corr( Z (s,\widetilde{\vk{v}}), Z (t, \widetilde{\vk{v}}'))\sim  K^2(|t-s|)+\frac{1}{4}\sum_{i=1}^{k-1}(v_i-v_i')^2
\EQN
as $ |t-s|\rw 0, |v_i-v_i'|\rw 0, 1\leq i \leq k-1.$ Let $a_0=(4^{-1}M_1)^{1/\alpha} $, and $a_1=a_2=\cdots=a_{k-1}=1/2$. Then, using the same arguments as in Lemma 2.1 in \cite{QueWat72} and in Lemma 2.2 in \cite{QuaWat73}, we conclude that
\BQNY
&&\pk{ \sup_{ t \in \widetilde{F}_{j,d} } \sup_{ \widetilde{\vk{v}} \in \widetilde{G}_{j,d}\cap [-\eta,\eta]^{k-1}}Z((4^{-1}M_1)^{1/\alpha} t,\widetilde{\vk{v}})> \sqrt{g(s_{j,d}^{(0)})}  }\\
&&\underset{\sim}\ge \frac{(2\eta)^{k-1} d }{q(g(s_{j,d}^{(0)})) (2h_{j,d} )^{k-1}} \frac{1}{\sqrt{2\pi} h_{j,d}} e^{-\frac{g(s_{j,d}^{(0)})}{2}}a_0 a_1^{k-1} \left(1-2^k \underset{ \vk n\neq \vk 0}{\sum_{0\le n_i<\IF, 1\le i\le k}} \left(1- \Phi\left( (n_1 a_0)^{\alpha/2}/2 + a_1/2\sum_{i=2}^k n_i \right)\right) \right)\\
&&=:  \mathcal{C }   \frac{(g(s_{j,d}^{(0)}) )^{k/2-1}}{q(g(s_{j,d}^{(0)}))} e^{-\frac{g(s_{j,d}^{(0)})}{2}} d
\EQNY
as $j\to\IF$, where $\Phi(\cdot)$ is the standard normal distribution function.
Therefore, for any $\delta, d$ small enough,
\BQN\label{con2}
\sum_{i=N_{d,\delta}}^\IF\mathbb{P}\left(W_{il_0,d}^c\right)
 &\geq& \mathcal{C }   \sum_{i=N_{d,\delta}}^\IF \frac{(g(s_{il_0,d}^{(0)}))^{\frac{k}{2}-1}}{ q(g(s_{il_0,d}^{(0)})) }e^{-\frac{g(s_{il_0,d}^{(0)})}{2}}d\nonumber\\
 &\geq& \mathcal{C }   \frac{1}{2l_0} \sum_{i=N_{d,\delta}}^{\IF}\int_{-(i+1)d}^{-(i+2)d}\frac{(g(\invf(t)))^{k/2-1}}{q((g(\invf(t))))}e^{-\frac{g(\invf(t))}{2}}dt\nonumber\\
 &\geq&  \mathcal{C }    \frac{1}{2l_0}\int_{0}^{\delta/2} (C(t))^{1/\alpha}\frac{(g(t))^{\frac{k}{2}-1}}{q(g(t)) }e^{-\frac{g(t)}{2}}dt=\IF.
\EQN
 Next, we prove
\BQN\label{Borel1}
\lim_{l\rw\IF}\left(\pk{\bigcap_{i=l}^\IF W_{il_0,d}}-\prod_{i=l}^\IF\pk{ W_{il_0,d}}\right)=0.
\EQN
 In view of Normal Comparison Lemma (see e.g., \cite{leadbetter1983extremes}), we have for any $m>l$
 \BQNY
 &&\left|\pk{\bigcap_{i=l}^m W_{il_0,d}}-\prod_{i=l}^m \pk{ W_{il_0,d}}\right|\\
 && \ \ \leq \frac{1}{2\pi}\sum_{l\leq i<j\le m}  \underset{{(t_1,\vk{v}_1)\in (F_{jl_0,d}\times (G_{jl_0,d}\cap \mathcal{S}_{k-1}))}}{\sum_{(t,\vk{v})\in (F_{il_0,d}\times (G_{il_0,d}\cap \mathcal{S}_{k-1}))}}  \frac{|r_{Y_k}(t,\vk{v}, t_1, \vk{v}_1)|}{\sqrt{1-r_{Y_k}^2(t,\vk{v}, t_1, \vk{v}_1)}}e^{-\frac{g^2(s_{il_0,d}^{(0)})+g^2(s_{jl_0,d}^{(0)})}{2(1+r_{Y_k}(t,\vk{v}, t_1, \vk{v}_1))}}=:\Lambda_{l,m}.
 \EQNY
Note that (\ref{B2}) implies that for any $0<\epsilon<\frac{\beta}{2(2-\beta)}$ with $\beta_1=\min (\beta, 1)$, there exist  $j_0 \in \mathbb{N}$ such that for any  $j> i \geq j_0$  and large enough $l_0$
 $$
\underset{{(t_1,\vk{v}_1)\in (F_{jl_0,d}\times (G_{jl_0,d}\cap \mathcal{S}_{k-1}))}}{\sup_{(t,\vk{v})\in (F_{il_0,d}\times (G_{il_0,d}\cap \mathcal{S}_{k-1}))}}  |r_{Y_k}(t,\vk{v}, t_1, \vk{v}_1))| \le \sup_{t\in \Delta^{(0)}_{il_0,d},t_1\in \Delta^{(0)}_{jl_0,d} }|r(t,t_1)| \leq M_0|(j-i)l_0|^{-\beta}\leq \epsilon.
 $$
With aid of (\ref{Restriction}), fundamental calculation yields that  for $l\geq j_0$,
\BQNY
\Lambda_{l,\IF}&\leq&\frac{M_0}{2\pi\sqrt{1-\epsilon^2}}\sum_{l\leq i< j<\IF}|(j-i)l_0|^{-\beta} \frac{d^2}{(q(g(s_{il_0,d}^{(0)}) ))^2}  \left(\sqrt{g(s_{il_0,d}^{(0)})}\right)^{2k} e^{-\frac{g(s_{il_0,d}^{(0)})+g(s_{jl_0,d}^{(0)})}{2(1+M_0|(j-i)l_0|^{-\beta})}}\\
 &\leq&  \mathcal{C }   \sum_{i=l}^\IF\left(g(s_{i,d}^{(0)})\right)^{k+4/\alpha}e^{-\frac{g(s_{i,d}^{(0)})}{2(1+\epsilon)}}\sum_{j=i+ l_0}^\IF|j-i|^{-\beta}\left(g(s_{j,d}^{(0)})\right)^{k+4/\alpha}e^{-\frac{g(s_{j,d}^{(0)})}{2(1+\epsilon)}}\\
 &\leq&  \mathcal{C }   \sum_{i=l}^\IF\left(\ln i\right)^{k+4/\alpha}i^{-\frac{1}{1+\epsilon}}\sum_{j=i+ l_0}^\IF|j-i|^{-\beta}\left(\ln j\right)^{k+4/\alpha}j^{-\frac{1}{1+\epsilon}}\\
 &\leq&   \mathcal{C }   \sum_{i=l}^\IF i^{-\frac{1}{1+2\epsilon}}\sum_{j=i+1}^\IF|j-i|^{-\beta}j^{-\frac{1}{1+2\epsilon}}\\
 &\leq&   \mathcal{C }   \sum_{i=l}^\IF i^{-\frac{1}{1+2\epsilon}}\left(\sum_{j=i+1}^{2i}|j-i|^{-\beta}j^{-\frac{1}{1+2\epsilon}}+\sum_{j=2i}^
 {\IF}|j-i|^{-\beta}j^{-\frac{1}{1+2\epsilon}}\right).
\EQNY
Since
$$\sum_{j=i+1}^{2i}|j-i|^{-\beta}j^{-\frac{1}{1+2\epsilon}}\leq \sum_{j=i+1}^{2i}|j-i|^{-\beta_1}j^{-\frac{1}{1+2\epsilon}}\leq  \mathcal{C }   i^{1-\beta_1-\frac{1}{1+2\epsilon}}$$
and
$$\sum_{j=2i}^
 {\IF}|j-i|^{-\beta}j^{-\frac{1}{1+2\epsilon}}\leq \sum_{j=2i}^
 {\IF}2^{\beta}j^{-\frac{1}{1+2\epsilon}-\beta}\leq  \mathcal{C }   i^{1-\frac{1}{1+2\epsilon}-\beta},$$
 then
\BQNY
\Lambda_{l,\IF}\leq  \mathcal{C }   \sum_{i=l}^\IF i^{1-\beta_1-\frac{2}{1+2\epsilon}}<\IF, \ \ l\geq l_0,
\EQNY
which implies that (\ref{Borel1}) holds. This completes the proof. \QED

\prooftheo{Addition}  An immediate application of Theorem \ref{THM} yields that 
$$
\pk{\sup_{t\in \mathcal{E}(S)}\left(\chi_{\vk{b}}^2(t)- g(t)\right)<\IF }=1
$$
provided that \aH{$|f(S)|<\IF$} or $I_g(S)<\IF$. %Moreover, by Theorem \ref{THM}, if $I_g(S)=\IF$, then
%$$\pk{\chi_{\vk{b}}^2(t)\leq g(t) \ \ \text{ultimately as}\ \ t\rw S }=0. $$
For the case $I_g(S)=\IF$, without loss of generality, we focus on $S=0$. For any such function $g(\cdot)$ satisfying $I_g(0)=\IF$, we can find a nonnegative continuous function $g_1(\cdot)$ (to be determined later) such that $g_1(t)\uparrow \IF$ as $t\rw 0$ and $I_{g+g_1}(0)=\IF$ hold. Then, by Theorem \ref{THM}  we have
\BQNY
\pk{\chi_{\vk{b}}^2(t)\leq g(t)+g_1(t) \ \mbox{ultimately as} \ t\rw 0}=0,
\EQNY
which implies that
\BQN  \label{01law_2}
\pk{\sup_{t\in (0,1/2]}\left(\chi_{\vk{b}}^2(t)- g(t)\right)=\IF }=1,
\EQN
and the proof will be complete.

Now we give one choice of the function $g_1(\cdot)$. Define $F(s)=\int_{s}^{1/2} (C(t))^{1/\alpha}\frac{(g(t))^{\frac{k}{2}-1}}{q(g(t)) }e^{-\frac{g(t)}{2}}dt$, and let $\invF(n)=\inf \{s\in (0,1/2]: F(s)=n\}, n\in \mathbb{N}$. Further, we construct a
nondecreasing function $w(\cdot)$ such that
$ w(t)=\frac{t-\invF(n)}{(n-1)\left(\invF(n-1)-\invF(n)\right)}+\frac{t-\invF(n-1)}{n\left(\invF(n)-\invF(n-1)\right)}, t\in [\invF(n), \invF(n-1)), n\geq 2$ and $w(t)=1, t\in [\invF(1), 1/2]$.
  Let $g_1(t)=\min(-2\ln w(t), g(t)), t\in(0,1/2]$. If $\alpha\in (0,2)$ or $\alpha=2, k>1$, then by the fact that
  $
  \frac{t^{\frac{k}{2}-1}}{q(t)}
  $
  is a regularly varying function at  $\IF$ with index $k/2-1+1/\alpha> 0$, we have
  \BQN\label{new}
  \frac{(g(t)+g_1(t))^{\frac{k}{2}-1}}{q( g(t)+g_1(t) )}\geq  \mathcal{C }  \frac{(g(t))^{\frac{k}{2}-1}}{q(g(t)) }
  \EQN
  for any $t\in (0,\overleftarrow{F}(N_0)]$ when $N_0$ is large enough. In light of {\bf B(0)}, one can easily check that (\ref{new}) also holds for $\alpha+2$ and $k=1$.
 Consequently,
\BQNY
&&\int_0^{1/2}(C(t))^{1/\alpha}\frac{(g(t)+g_1(t))^{\frac{k}{2}-1}}{q  (g(t)+g_1(t)) }e^{-\frac{g(t)+g_1(t)}{2}}dt\geq  \mathcal{C }  \int_0^{\overleftarrow{F}(N_0)}(C(t))^{1/\alpha}\frac{(g(t))^{\frac{k}{2}-1}}{q(g(t)) }e^{-\frac{g(t)}{2}}w(t)dt\\
&&\ \ =  \mathcal{C }  \sum_{n=N_0+1}^\IF\int_{\invF(n)}^{\invF(n-1)}(C(t))^{1/\alpha}\frac{(g(t))^{\frac{k}{2}-1}}{q(g(t)) }e^{-\frac{g(t)}{2}}w(t)dt
\geq   \mathcal{C }  \sum_{n=N_0+1}^{\IF} \frac{1}{n}=\IF,
\EQNY
which completes the proof. \QED

\bigskip

{\bf Acknowledgement}: We are grateful to the referees for their comments and suggestions which significantly improved the manuscript.
 P. Liu acknowledges partial
support from the Swiss National Science Foundation Project 200021-166274 and L. Ji acknowledges partial
support from the Swiss National Science Foundation Project P300P2-167649.
%and the project RARE -318984  (an FP7  Marie Curie IRSES Fellowship).
\bibliographystyle{plain}

 \bibliography{GCCD}

\end{document}